\documentclass[oneside,notitlepage,12pt]{article}

\usepackage{amssymb}
\usepackage[leqno]{amsmath}
\usepackage{amsfonts}
\usepackage{amsopn}
\usepackage{amstext}
\usepackage{amsthm}

\usepackage[all]{xy}

\usepackage[colorlinks,backref]{hyperref}
\usepackage{makeidx}

\usepackage{calrsfs}
\usepackage{fourier-orns}
\usepackage{hieroglf} 
\usepackage{clock} \ClockFrametrue\ClockStyle1

\newenvironment{pf}{\begin{proof}}{\end{proof}}

\newcommand{\Ef}{{\cal{F}}}

\newcommand{\Nat}{{\mathbb{N}}}
\newcommand{\Qyu}{{\mathbb{Q}}}
\newcommand{\Err}{{\mathbb{R}}}

\newcommand{\U}{{\mathbb U}}
\newcommand{\G}{{\mathbb G}}

\newcommand{\sig}{\sigma}
\newcommand{\eps}{\varepsilon}
\renewcommand{\phi}{\varphi}
\renewcommand{\rho}{\varrho}

\newcommand{\rest}{\restriction}
\newcommand{\unii}{\mathbb{I}}
\newcommand{\ntr}{n\in\omega}
\newcommand{\Ntr}{n\in{\mathbb{N}}}
\newcommand{\loe}{\leqslant}
\newcommand{\goe}{\geqslant}

\newcommand{\subs}{\subseteq}

\newcommand{\nnempty}{\ne\emptyset}

\newcommand{\ovr}{\overline}

\newcommand{\dist}{\operatorname{dist}}

\newcommand{\id}[1]{{\operatorname{i\!d}_{#1}}} 

\newcommand{\dom}{\operatorname{dom}}
\newcommand{\cod}{\operatorname{cod}}

\newcommand{\liminv}{\varprojlim}

\newcommand{\oraz}{\qquad\text{and}\qquad}

\newcommand{\Lip}[1]{\operatorname{Lip}\left(#1\right)}

\newcommand{\poset}{{\mathbb{P}}}

\newcommand{\ob}[1]{\operatorname{Ob}\left(#1\right)}

\newcommand{\Metrics}{\mathfrak M\!\mathfrak s}
\newcommand{\Komp}{\mathfrak C \!\mathfrak o \!\mathfrak m \!\mathfrak p}

\newcommand{\by}[1]{/_{#1}}

\newtheorem{tw}{Theorem}[section]
\newtheorem{wn}[tw]{Corollary}
\newtheorem{lm}[tw]{Lemma}
\newtheorem{prop}[tw]{Proposition}

\theoremstyle{definition}

\newtheorem{ex}[tw]{Example}

\theoremstyle{remark}
\newtheorem{uwgi}[tw]{Remark}

\providecommand{\nat}{\omega}

\newcommand{\setof}[2]{\{#1\colon #2\}}

\newcommand{\sett}[2]{\{#1\}_{#2}}
\newcommand{\bigsett}[2]{\left\{#1\right\}_{#2}}
\newcommand{\sn}[1]{\{#1\}} 
\newcommand{\dn}[2]{\{#1,#2\}} 
\newcommand{\pair}[2]{{\langle{#1},{#2}\rangle}} 
\newcommand{\triple}[3]{\langle #1, #2, #3 \rangle} 
\newcommand{\map}[3]{#1\colon #2 \to #3} 
\newcommand{\img}[2]{#1[#2]} 
\newcommand{\inv}[2]{{#1}^{-1}[#2]} 

\newcommand{\ciag}[1]{{\sett{{#1}_n}{\ntr}}}

\newcommand{\iso}{\approx}

\newcommand{\anorm}{\|\cdot\|}
\newcommand{\norm}[1]{\|#1\|}

\newcommand{\fK}{{\mathfrak{K}}}

\newcommand{\fB}{{\mathfrak{B}}}
\newcommand{\fC}{{\mathfrak{C}}}

\newcommand{\fI}{{\mathfrak{I}}}

\newcommand{\fra}{Fra\"iss\'e}

\newcommand{\bS}{{\mathbb{S}}}

\newcommand{\cmp}{\circ} 

\newcommand{\wek}[1]{{\vec{#1}}}

\newcommand{\lcat}[1]{{#1}^{\heartsuit}}
\newcommand{\rcat}[1]{{#1}}
\newcommand{\paircats}[1]{\pair {\lcat {#1}}{\rcat {#1}}}

\newcommand{\Cantor}{2^\omega}

\newcommand{\separator}{\begin{center} \leafright \leafright \decotwo \leafleft \leafleft \end{center}}

\newcommand{\define}[2]{{\em #1}\index{#2}}

\newcommand{\ciagi}[1]{\sig{#1}}
\newcommand{\ciagilr}[1]{{\sig({\lcat #1},{\rcat #1})}}

\newcommand{\Gurarii}{Gurari\u\i}

\newcommand{\Ban}{\mathfrak B^{\operatorname{fd}}}

\title{Metric-enriched categories and approximate \fra\ limits}

\author{
Wies{\l}aw Kubi\'s
\footnote{Research supported by the NCN grant DEC-2011/03/B/ST1/00419.}
\\
{\small Institute of Mathematics, Jan Kochanowski University in Kielce, Poland}\\
{\small Institute of Mathematics, Academy of Sciences of the Czech Republic}\\
{\small (kubis@math.cas.cz)}
}

\date{\clocktime\ \today}

\makeindex

\begin{document}

\maketitle

\begin{abstract}
We develop the theory of \define{approximate}{approximate \fra\ limit} \fra\ limits in the context of categories enriched over metric spaces.
Among applications, we construct a generic projection on the \Gurarii\ space and we present a simpler proof of a recent characterization of the pseudo-arc, due to Irwin and Solecki.


\end{abstract}

\tableofcontents

\section{Introduction}

Suppose that we are given a category $\fK$ of some sort of embeddings; let us say that the objects of $\fK$ are \emph{small}.
Now, assume that certain sequences in $\fK$ (i.e. covariant functors from the set of natural numbers into $\fK$) have co-limits in a bigger category, denoted by $\ovr \fK$; the objects of $\ovr \fK$ will be called \emph{big}.
We would like to have a ``generic" sequence in $\fK$, whose co-limit (the \emph{$\ovr \fK$-generic object}) will accommodate all $\fK$-objects, and which will have the best possible homogeneity property.
Specifically, a $\ovr \fK$-object $U$ is \emph{$\fK$-homogeneous} if given a $\fK$-object $a$, given $\fK$-arrows $\map {e_0,e_1} a U$, there exists an automorphism $\map h U U$ such that $h \cmp e_0 = e_1$.

Unfortunately, it turns out that there exists a very natural category $\Ban$, namely, all finite-dimensional Banach spaces with linear isometric embeddings, where the generic object is the \Gurarii\ space which is only \emph{almost} homogeneous with respect to $\Ban$.
A general back-and-forth type argument says that the generic object is unique, therefore no generic separable Banach space can be $\Ban$-homogeneous.

Generic objects are in fact straight generalizations of \emph{\fra\ limits}, well known in model theory.
The example mentioned above would fit into the \fra\ theory, however there are continuum many isometric types of finite-dimensional Banach spaces, therefore the existence of a generic sequence is not so obvious.

Our aim is to find a general framework for these situations.
It turns out that dealing with a category enriched over metric spaces is a possible solution.
On the other hand, referring to the example of Banach spaces, one quickly realizes that the \Gurarii\ space cannot be fully explained in the category of isometric embeddings.
The reason is precisely the fact that it lacks the homogeneity property: There exist two 1-dimensional subspaces for which no bijective linear isometry maps the first one onto the other.
Thus, it is necessary to use isomorphic embeddings or, in other words, to work in the category of non-expansive linear embeddings between finite-dimensional spaces.

The purpose of this work is developing the theory of \emph{approximate \fra\ limits} in the context of categories enriched over metric spaces.
One of the features is that we deal with such a category $\rcat \fK$ together with a (usually, very natural) subcategory $\lcat \fK$ with the same objects as $\rcat \fK$, such that all $\lcat \fK$-sequences have co-limits in $\ovr \fK$, and $\lcat \fK$ satisfies some natural assumptions good enough for the existence of a generic sequence.
This sequence leads to the \fra\ limit $U$ in $\ovr \fK$, which is \emph{almost $\fK$-homogeneous} in the sense that given a $\fK$-object $a$, given $\ovr {\lcat \fK}$-arrows $\map {e_0,e_1} a U$, given $\eps > 0$, there exists an automorphism $\map h U U$ in $\ovr {\lcat \fK}$ such that $h \cmp e_0$ is $\eps$-close to $e_1$.
As mentioned above, all these categories are enriched over metric spaces which means, roughly speaking, that a distance $\rho(f,g)$ is defined for each pair of arrows $f,g$ with common domain and common co-domain, allowing us to measure the ``commutativity" of diagrams.

A prototype example is the category $\rcat \Metrics$ of metric spaces with $1$-Lipschitz mappings, where $\lcat \Metrics$ consists of all isometric embeddings.
Actually, this is the canonical example of a category enriched over itself.
Another example is the category $\rcat \Komp$ of $1$-Lipschitz maps between metric compacta, where the arrows are formally reversed: an arrow from $X$ to $Y$ is a $1$-Lipschitz map from $Y$ into $X$.
This is just because we would like to consider inverse sequences whose limits, in the category of topological spaces, are again compact metrizable spaces.
Here, $\lcat\Komp$ will be the subcategory of all quotient maps in $\rcat \Komp$.

One has to admit that every category $\fK$ is enriched over metric spaces by the zero-one metric.
In this case, our main results are actually part of the category-theoretic \fra\ theory developed in \cite{DrGoe92} and \cite{KubFrais}.

One has to mention that a parallel research in approximate \fra\ theory has been recently done by Ben Yaacov~\cite{BenY} in continuous model theory (partially inspired by the PhD thesis of Schoretsanitis~\cite{Schor}).
Ben Yaacov's main tool is the concept of \emph{bi-Kat\v etov maps}, smartly encoding almost isometric embeddings of metric structures.

The concept of categories enriched over metric spaces goes back to Eilenberg \& Steenrod~\cite{EilSte}, although one of the main inspirations for our study comes from a note of Mioduszewski \cite{Miod1} on $\eps$-commuting diagrams and inverse limits of compact metric spaces.

The paper is organized as follows.
Section~\ref{Sectwotuu} contains the main definitions, like metric-enriched and norm category, sequences, almost amalgamation property, etc.
Section~\ref{SecGenerics} contains the main results, starting with the crucial notion of a \fra\ sequence, characterizing its existence and proving its main properties.
Section~\ref{SectAppsy} contains selected applications: a simple description of the \Gurarii\ space, universal projections, a new point of view on the Cantor set, and finally a discussion of the pseudo-arc, including its new characterization.

\section{Metric-enriched categories}\label{Sectwotuu}

Let $\Metrics$ denote the category of metric spaces with non-expansive (i.e., $1$-Lipschitz) mappings.
A category $\fK$ is \define{enriched over}{enriched category} $\Metrics$ if for every $\fK$-objects $a, b$ there is a metric $\rho$ on the set of $\fK$-arrows $\fK(a,b)$ so that the composition operator is non-expansive on both sides.
More precisely, we have
\begin{equation}
\rho(f_0 \cmp g, f_1 \cmp g) \loe \rho(f_0, f_1) \oraz \rho(h \cmp f_0, h \cmp f_1) \loe \rho(f_0, f_1)
\tag{M}\label{Equametric}
\end{equation}\define{}{(M)}
whenever the compositions make sense.
This allows us to consider $\eps$-commutative diagrams, with the obvious meaning.
Formally, on each hom-set we have a different metric, but there is no ambiguity with using always the same letter $\rho$.
For the sake of convenience, we allow $+\infty$ as a possible value of the metric.

Later on, we shall say that $\fK$ is \define{metric-enriched}{metric-enriched category} having in mind that $\fK$ is enriched over $\Metrics$.

\subsection{Norms}

From now on, we fix a pair $\paircats \fK$, where $\rcat \fK$ is a metric-enriched category and $\lcat \fK$ is its subcategory with the same objects.

Given $f \in \rcat \fK$, define
$$\mu(f) = \inf \setof{ \rho(j \cmp f, i) }{i, j \in \lcat \fK},$$
where only compatible arrows $i, j$ are taken into account.
We allow the possibility that $\mu(f) = +\infty$, in particular when there are no $i,j \in \lcat \fK$ with $\dom(j) = \cod(f) = \cod(i)$ and $\dom(i) = \dom(f)$.
Obviously,
$$\lcat \fK \subs \setof{f\in \rcat \fK}{\mu(f) = 0},$$
and in typical cases the equality holds.
The meaning of $\mu(f)$ is ``measure of distortion".
We call $\mu$ the \define{norm}{norm} induced by $\triple {\lcat \fK}{\rcat \fK}\rho$.
An arrow $f$ satisfying $\mu(f) = 0$ will be called a \define{0-arrow}{0-arrow}.
It turns out that the inverse of a 0-isomorphism is a 0-arrow:

\begin{prop}
Assume $\map h a b$ is an isomorphism. Then $\mu(h) = \mu(h^{-1})$.
\end{prop}

\begin{pf}
Given $\lcat \fK$-arrows $\map i a c$, $\map j b c$, we have
$$\rho(j \cmp h, i) = \rho(j \cmp h, i \cmp h^{-1} \cmp h) \loe \rho(j, i \cmp h^{-1}).$$
By symmetry, $\rho(j \cmp h, i) = \rho(j, i \cmp h^{-1})$.
Taking the infimum for all possible compatible $\lcat \fK$-arrows $i,j$ we obtain $\mu(h) = \mu(h^{-1})$.
\end{pf}

From now on, the triple $\triple {\lcat \fK}{\rcat \fK}\rho$ will be called a \define{normed category}{normed category}.
We shall usually omit the metric $\rho$, just saying that $\paircats \fK$ is a normed category.
It is clear that the norm (always denoted by $\mu$) is determined by the pair $\triple {\lcat \fK} {\rcat \fK}\rho$.

\begin{ex}\label{Exebrigo}
Let $\rcat \fB$ be the category of all finite-dimensional Banach spaces with linear operators of norm $\loe 1$ and let $\lcat \fB$ be the subcategory of all isometric embeddings.
It is clear that $\rcat \fB$ is a metric-enriched category, where $\rho(f,g) =  \norm{f - g}$ for $\map {f,g} X Y$ in $\rcat \fB$.

Now assume $\map f X Y$ is a linear operator satisfying
\begin{equation}
(1 - \eps) \norm x \loe \norm{ f(x) } \loe \norm x.
\tag{$*$}\label{Eqtjojhre}
\end{equation}
We claim that $\mu(f) \loe \eps$.
In fact, consider $Z = X \oplus Y$ with the norm defined by the following formula:
$$\norm{\pair x y} = \inf \setof{ \norm{u}_X + \norm{v}_Y + \eps \norm{w}_X }{ \pair x y = \pair u v + \pair w {-f(w)}, \; u,w \in X,\; v \in Y }.$$
Let $\map i X Z$, $\map j Y Z$ be the canonical injections.
Note that $\norm {i - j \cmp f} \loe \eps$, just by definition.
It remains to check that $i$, $j$ are isometric embeddings, that is, they belong to the category $\lcat \fB$.

It is clear that $\norm i \loe 1$ and $\norm j \loe 1$.
On the other hand, if $x = u + w$ and $0 = v - f(w)$ then
$$\norm{u}_X + \norm{v}_Y + \eps \norm{w}_X \goe \norm{u}_X + (1-\eps)\norm{w}_X + \eps\norm{w}_X \goe \norm{u + w}_X = \norm{x}_X.$$
This shows that $\norm{i(x)} = \norm{x}_X$.
A similar calculation shows that $\norm{j(y)} = \norm{y}_Y$.

Thus, if $f$ satisfies (\ref{Eqtjojhre}) then $\mu(f)\loe \eps$.
On the other hand, if $\norm {f(x)} = 1 - \eps$ for some $x$ with $\norm x = 1$, then given isometric embeddings $i$, $j$, we have
$$\norm {i(x) - j(f(x))} \goe | \norm {i(x)} - \norm{j(f(x))} | = |1 - \norm {f(x)}| = \eps.$$
Finally, we conclude that $\mu(f) = \eps$, where $\eps \goe 0$ is minimal for which the inequality (\ref{Eqtjojhre}) holds.
\end{ex}

The example above is taken from \cite{GarKub0}, where it is proved that the embeddings $i,j$ lead to a universal object in the appropriate category.
A slightly more technical argument is given in \cite{KubSol}.

\subsection{Sequences and approximate arrows}

In order to speak about ``big" objects, we shall introduce a natural category of sequences.
By this way, ``big" objects will be identified with (equivalence classes of) sequences, having in mind their co-limits.
For example, every separable complete metric space is the completion of the union of a chain of finite metric spaces, therefore it can be described in terms of sequences in the category $\Metrics$.
In general, we consider sequences in which the bonding arrows are some sort of monics, therefore it is convenient to work in a pair of categories.
For our applications, it is sufficient to discuss the case of normed categories.

Let $\paircats \fK$ be a fixed normed category.
A \define{sequence}{sequence} in $\fK$ is formally a covariant functor from the set of natural numbers $\nat$ into $\fK$.

Denote by $\ciagilr \fK$ the category of all sequences in $\lcat \fK$ with arrows in $\rcat \fK$.
This is indeed a category with arrows being equivalence classes of semi-natural transformations.
A \define{semi-natural transformation}{semi-natural transformation} from $\wek x$ to $\wek y$ is, by definition, a natural transformation from $\wek x$ to $\wek y \cmp \phi$ for some increasing function $\map \phi \nat \nat$.
Slightly abusing notation, we shall consider transformations (arrows) of sequences, having in mind their equivalence classes.
Thus, an arrow from $\wek x$ to $\wek y$ is a sequence of arrows $\wek f = \ciag f$ together with an increasing map $\map \phi \nat \nat$ such that for each $n < m$ the diagram
$$\xymatrix{
y_{\phi(n)} \ar[rr]^{y_{\phi(n)}^{\phi(m)}} & & y_{\phi(m)} \\
x_n \ar[u]^{f_n} \ar[rr]_{x_n^m} & & x_m \ar[u]_{f_m}
}$$
is commutative.

Since the category $\rcat \fK$ is enriched over $\Metrics$, it is natural to allow more arrows in $\ciagilr \fK$.
An \define{approximate arrow}{approximate arrow} from a sequence $\wek x$ into a sequence $\wek y$ is $\wek f = \ciag f \subs \rcat \fK$ together with an increasing map $\map \phi \nat \nat$ satisfying the following condition: 
\begin{enumerate}
\item[($\circlearrowleft$)] For every $\eps > 0$ there exists $n_0$ such that all diagrams of the form
$$\xymatrix{
y_{\phi(n)} \ar[rr]^{y_{\phi(n)}^{\phi(m)}} & & y_{\phi(m)} \\
x_n \ar[u]^{f_n} \ar[rr]_{x_n^m} & & x_m \ar[u]_{f_m}
}$$
\define{}{($\circlearrowleft$)}
are $\eps$-commutative, i.e., $\rho(f_m \cmp x_n^m, y_{\phi(n)}^{\phi(m)} \cmp f_n) < \eps$ whenever $n_0 \loe n < m$.
\end{enumerate}
It is obvious that the composition of approximate arrows is an approximate arrow, therefore $\ciagilr \fK$ is indeed a category.
It turns out that $\ciagilr \fK$ is naturally metric-enriched.
Indeed, given approximate arrows $\map {\wek f}{\wek x}{\wek y}$, $\map {\wek g}{\wek x}{\wek y}$, define
\begin{equation}
\rho(\wek f, \wek g) := \lim_{n\to\infty}\lim_{m > n} \rho(y_n^m \cmp f_n, y_n^m \cmp g_n),
\tag{$\rightrightarrows$}
\end{equation}\define{}{($\rightrightarrows$)}
assuming that $\wek f$ and $\wek g$ are natural transformations.
Replacing $\wek y$ by its cofinal subsequence, we can make such assumption, without loss of generality.
We need to show that the limit above exists.

Given $n < m < \ell$, we have
$$\rho(y_n^\ell \cmp f_n, y_n^\ell \cmp g_n) = \rho(y_m^\ell \cmp y_n^m \cmp f_n, y_m^\ell \cmp y_n^m \cmp g_n) \loe \rho(y_n^m \cmp f_n, y_n^m \cmp g_n),$$
therefore for each $\ntr$ the sequence $\sett{\rho(y^m_n \cmp f_n, y^m_n \cmp g_n)}{m>n}$ is decreasing.
On the other hand, given $\eps>0$, given $n_0 \loe n < k < m$ such that ($\circlearrowleft$) holds for both $\sett{f_n}{n \goe n_0}$ and $\sett{g_n}{n \goe n_0}$, we have that
\begin{align*}
\rho(y^m_n \cmp f_n, y^m_n \cmp g_n) &\loe \rho(y^m_k \cmp y^k_n \cmp f_n, y^m_k \cmp f_k \cmp x^k_n) + \rho(y^m_k \cmp f_k \cmp x^k_n, y^m_k \cmp g_k \cmp x^k_n) \\
&+ \rho( y^m_k \cmp g_k \cmp x^k_n, y^m_k \cmp y^k_n \cmp g_n)
\loe \rho(y^m_k \cmp f_k, y^m_k \cmp g_k) \\
&+ \rho(y^k_n \cmp f_n, f_k \cmp x^k_n) + \rho(g_k \cmp x^k_n, y^k_n \cmp g_n) \\
&< \rho(y^m_k \cmp f_k, y^m_k \cmp g_k) + 2\eps.
\end{align*}
Passing to the limit as $m \to \infty$, we see that the sequence $$\bigsett{\lim_{m > n} \rho(y_n^m \cmp f_n, y_n^m \cmp g_n)}{\ntr}$$ is increasing.
This shows that the double limit in ($\rightrightarrows$) exists.

One should mention that a much more natural definition for $\rho$ would be
\begin{equation}
\rho(\wek f, \wek g) = \lim_{n \to \infty}\rho(f_n, g_n).
\tag{$\twoheadrightarrow$}
\end{equation}\define{}{($\twoheadrightarrow$)}
The problem is that this limit may not exist in general.
In practice however, we shall always have
$$\rho(i \cmp f, i \cmp g) = \rho(f,g)$$ whenever $i \in \lcat \fK$.
To be more precise, when $\fK$ is a metric-enriched category, an arrow $i \in \fK$ is called a \define{monic}{monic} (or a \emph{monomorphism}) if the above equation holds for arbitrary compatible arrows $f,g \in \fK$.
This is an obvious generalization of the notion of a monic in category theory.
In fact, every category is metric-enriched over the 0-1 metric: $\rho(f,g) = 0$ iff $f = g$.

In all natural examples of normed categories $\paircats \fK$ the subcategory $\lcat \fK$ consists of arrows that are monics in $\rcat \fK$.
Thus, we can use ($\twoheadrightarrow$) instead of the less natural ($\rightrightarrows$) as the definition of $\rho(\wek f, \wek g)$.
It is an easy exercise to check that $\rho$ is indeed a metric on each hom-set of $\ciagilr \fK$ and that all composition operators are non-expansive.

The important fact is that every normed category $\paircats \fK$ naturally embeds into $\ciagilr \fK$, identifying a $\lcat \fK$-object $x$ with the sequence of identities
$$\xymatrix{
x \ar[r] & x \ar[r] & x \ar[r] & \cdots
}$$
and every $\rcat \fK$-arrow becomes a natural transformation between sequences of identities.
Actually, it may happen that $\rcat \fK$ is not a full subcategory of $\ciagilr \fK$.
Namely, every Cauchy sequence of $\rcat \fK$-arrows $\map {f_n}x y$ is an approximate arrow from $x$ to $y$, regarded as sequences.
Now, if $\rcat \fK(x,y)$ is not complete, there may be no $\map f x y$ satisfying $f = \lim_{n \to \infty} f_n$.
The problem disappears if all hom-sets of $\rcat \fK$ are complete with respect to the metric $\rho$.

\separator

In practice, we can partially ignore the construction described above, because usually there is a canonical faithful functor from $\ciagilr \fK$ into a natural category containing $\rcat \fK$ and in which all countable sequences in $\lcat \fK$ have co-limits.
Two relevant examples are described below.

\begin{ex}
Let $\rcat \fK$ be the category of all finite metric spaces with non-expansive mappings and let $\lcat \fK$ be the subcategory of isometric embeddings.
Let $\fC$ be the category of all separable complete metric spaces.
Given a $\lcat \fK$-sequence $\wek x$, we can identify it with a chain of finite metric spaces, therefore it is natural to consider $\lim \wek x$ to be the completion of the union of this chain.
This is in fact the co-limit of $\wek x$ in the category $\Metrics$.
In particular, $\lim \wek x = \lim \wek y$, whenever $\wek x$ and $\wek y$ are equivalent.
Furthermore, every approximate arrow $\map {\wek f} {\wek x}{\wek y}$ ``converges" to a non-expansive map $\map F {\lim \wek x}{\lim \wek y}$ and again it is defined up to an equivalence of approximate arrows.
By this way we have defined a canonical faithful functor $\map \lim {\ciagilr \fK} \fC$.
Unfortunately, since we are restricted to 1-Lipschitz mappings, the functor $F$ is not onto.
The simplest example is as follows.

Let $X = \dn 01 \cup \setof{\pm 1/n}{n\in\Nat}$ with the metric induced from the real line.
Let $X_n = \setof{\pm 1/k}{k < n}$ and let $Y = \dn 01$.
Then $X = \lim_{\ntr} X_n$ although the canonical embedding $\map e YX$ is not the co-limit of any approximate arrow from the sequence $\ciag X$ into $Y$ (the space $Y$ can be treated as the infinite constant sequence with identities).

Note that if we change $\fC$ to the category of all countable metric spaces then 
the canonical co-limit is just the union of the sequence, however
some approximate arrows of sequences would not have co-limits.
\end{ex}

Concerning applications, the situation where the canonical ``co-limiting" functor is not surjective does not cause any problems, since our main results say about the existence of certain arrows (or isomorphisms) of sequences only.
In the next example we have a better situation.

\begin{ex}\label{Extrgergx}
Let $\rcat \fB$ be, as in Example~\ref{Exebrigo}, the category of all finite-dimensional Banach spaces with linear operators of norm $\loe 1$ and let $\lcat \fB$ be the subcategory of all isometric embeddings.

Again, we have a canonical functor $\map \lim {\ciagilr \fB}\fC$, where $\fC$ is the category of all separable Banach spaces with non-expansive linear operators.
It turns out that this functor is surjective.
Namely, fix two Banach spaces $X = \ovr{\bigcup_{\ntr}X_n}$ and $Y = \ovr{\bigcup_{\ntr}Y_n}$, where $\ciag X$ and $\ciag Y$ are chains of finite-dimensional spaces.
Fix a linear operator $\map T X Y$ such that $\norm T \loe 1$ and let $T_n = T \rest X_n$.
By an easy induction, we define a sequence of linear operators $\map {T_n'}{X_n}{Y_n}$ so that $T_{n+1}'$ extends $T_n'$ and $\norm {T_n - T_n'} < 1/n$ for each $\ntr$.
Note that $\norm{T_n'} \loe 1 + 1/n$. Define $T''_n = \frac n{n+1}T_n'$.
Now $\wek t = \sett{T''_n}{\ntr}$ is a sequence of $\rcat \fB$-arrows and standard calculation shows that $\norm {T''_n - T_n} < 2/n$ for every $\ntr$.
Thus $\wek t$ is an approximate arrow from $\ciag X$ to $\ciag Y$ with $\lim \wek t = T$. 
\end{ex}

It is natural to extend the norm $\mu$ to the category of sequences.
More precisely, given an approximate arrow $\map {\wek f}{\wek x}{\wek y}$, we define
$$\mu(\wek f) = \lim_{n \to \infty}\mu(f_n).$$
This is indeed well defined, because given $\eps > 0$ and taking $n < m$ as in ($\circlearrowleft$), we have that $\mu(f_n) \loe \mu(f_m) + \eps$.
The function $\mu$ obviously extends the norm of $\paircats \fK$, although it is formally not a norm, because it is defined in a different way, without referring to any subcategory of $\ciagilr \fK$.
An approximate arrow $\wek f$ is a \define{0-arrow}{0-arrow of sequences} if $\mu(\wek f) = 0$.
We shall be interested in 0-arrows only.

\subsection{Almost amalgamations}

Let $\fK$ be a metric-enriched category.
We say that $\fK$ has \define{almost amalgamation property}{almost amalgamation property} if for every $\fK$-arrows $\map f c a$, $\map g c b$, for every $\eps > 0$, there exist $\fK$-arrows $\map {f'} a w$, $\map {g'} b w$ such that the diagram
$$\xymatrix{
b \ar[r]^{g'} & w \\
c \ar[u]^g \ar[r]_f & a \ar[u]_{f'}
}$$
is $\eps$-commutative, i.e. $\rho(f' \cmp f, g' \cmp g) < \eps$.
We say that $\fK$ has the \define{strict amalgamation property}{amalgamation property!-- strict} if for each $f,g$ the diagram above is commutative (i.e. no $\eps$ is needed).

It turns out that almost amalgamations can be moved to the bigger category $\rcat \fK$. Namely:

\begin{prop}\label{Pamelias}
Let $\paircats \fK$ be with the almost amalgamation property.
Then for every $\eps, \delta > 0$, for every $\rcat \fK$-arrows $\map f c a$, $\map g c b$ with $\mu(f) < \eps$, $\mu(g) < \delta$, there exist $\lcat \fK$-arrows $\map {f'} a w$ and $\map {g'} b w$ such that the diagram
$$\xymatrix{
b \ar[r]^{g'} & w \\
c \ar[r]_f \ar[u]^g & a \ar[u]_{f'}
}$$
is $(\eps + \delta)$-commutative.
\end{prop}

\begin{pf}
Fix $\eta > 0$.
Find $\lcat \fK$-arrows $i$, $j$ such that $\rho(j \cmp f, i) < \mu(f) + \eta$.
Find $\lcat \fK$-arrows $k$, $\ell$ such that $\rho(\ell \cmp g, k) < \mu(g) + \eta$.
Using the almost amalgamation property, find $\lcat \fK$-arrows $j'$, $\ell'$ such that $\rho(j' \cmp i, \ell' \cmp k) < \eta$.
Define $f' := j' \cmp j$, $g' := \ell' \cmp \ell$.
Then
\begin{align*}
\rho(f' \cmp f, g' \cmp g) &\loe \rho(j' \cmp j \cmp f, j' \cmp i) + \rho(j' \cmp i, \ell' \cmp k) + \rho(\ell' \cmp k, \ell' \cmp \ell \cmp g) \\
&< \rho(j \cmp f, i) + \eta + \rho(k, \ell \cmp g) < \mu(f) + \mu(g) + 3\eta.
\end{align*}
Thus, it is clear that if $\eta$ is small enough then $\rho(f' \cmp f, g' \cmp g) < \eps + \delta$.
\end{pf}

One more property needed later is that every two objects can be ``mapped" into a common one.
Namely, we say that a category $\fK$ is \define{directed}{directed category} if for every $a, b \in \ob \fK$ there is $c \in \ob \fK$ such that both hom-sets $\fK(a,c)$ and $\fK(b,c)$ are nonempty.
In model theory, this is usually called the \define{joint embedding property}{joint embedding property}.

In the context of normed categories, we are interested in directedness of the smaller category.
In the presence of amalgamations, this turns out to be equivalent to the directedness of the bigger category.

Let us say that a category $\fK$ has the \define{pseudo-amalgamation property}{pseudo-amalgamation} if for every $\fK$-arrows $\map f c a$, $\map g c b$ there exists $d \in \ob \fK$ such that both hom-sets $\fK(a,d)$, $\fK(b,d)$ are nonempty.

\begin{prop}
Let $\paircats \fK$ be a normed category such that $\lcat \fK$ has the pseudo-amalgamation property.
The following statements are equivalent:
\begin{enumerate}
	\item[{\rm(a)}] $\lcat \fK$ is directed.
	\item[{\rm(b)}] Given $a ,b \in \ob {\lcat \fK}$, there exist $\rcat \fK$-arrows $\map f a d$, $\map g b d$ such that $\mu(f) < +\infty$ and $\mu(g) < +\infty$.
\end{enumerate}
\end{prop}

\begin{pf}
Evidently, (a)$\implies$(b), because $\lcat \fK$ and $\rcat \fK$ have the same objects.
Suppose (b) holds and fix $a,b \in \ob {\lcat \fK}$.
Fix $f$, $g$ as in (b).
Since $\mu(f) < +\infty$, there exist $\lcat \fK$-arrows $\map i a v$, $\map j d v$ such that $\rho(i, j \cmp f) < +\infty$.
Similarly, there exist $\lcat \fK$-arrows $\map k b w$, $\map \ell d w$ such that $\rho(k, \ell \cmp g) < +\infty$.
Finally, using the pseudo-amalgamation property, find $z \in \ob {\lcat \fK}$ and some $\lcat \fK$-arrows $\map {j'} v z$ and $\map {\ell'} w z$.
The $\lcat \fK$-arrows $\map {j' \cmp i} a z$ and $\map {\ell' \cmp k} b z$ witness the directedness of $\lcat \fK$.
\end{pf}

\begin{prop}
Assume $\paircats \fK$ is a normed category with the almost amalgamation property. Then for every compatible $\rcat \fK$-arrows $f,g$ the following inequalities hold:
\begin{enumerate}
	\item[$(N_1)$] $\mu(f \cmp g) \loe \mu(f) + \mu(g)$.
\define{}{($N_1$)}
	\item[$(N_2)$] $\mu(g) \loe \mu(f) + \mu(f \cmp g)$.
\define{}{($N_2$)}
\end{enumerate}
\end{prop}

\begin{pf}
Fix $\eps > 0$ and fix $i,j,k,\ell \in \lcat \fK$ such that
$$\rho(j \cmp f, i) < \mu(f) + \eps/3 \oraz \rho(\ell \cmp g, k) < \mu(g) + \eps/3.$$
Using the almost amalgamation property, we can find $i', \ell' \in \lcat \fK$ such that
$$\rho(i' \cmp i, \ell' \cmp \ell) < \eps/3.$$
Combining these inequalities, we obtain that $\rho(\ell' \cmp k, i' \cmp j) < \mu(f) + \mu(g) + \eps$, which shows ($N_1$).

A similar argument shows ($N_2$).
\end{pf}

It is well-known and easy to see that the category of finite metric spaces with isometric embeddings has the strict amalgamation property.
The same result, whose precise formulation is given below, holds for Banach spaces. Most likely it belongs to the folklore, although we refer to~\cite{ACCGM} for a discussion of a more general statement.

\begin{prop}\label{PrpPuszautBanachs}
Let $\map i Z X$, $\map j Z Y$ be linear isometric embeddings of Banach spaces.
Then there exist linear isometric embeddings $\map {i'} X W$, $\map {j'} Y W$ such that
$$\xymatrix{
Y \ar[rr]^{j'} & & W \\
Z \ar[u]^{j} \ar[rr]_{i} & & X \ar[u]_{i'}
}$$
is a pushout square in the category of Banach spaces with linear operators of norm $\loe 1$.
Furthermore, $W = (X \oplus Y) \by \Delta$, where $\Delta = \setof{\pair {i(z)}{-j(z)}}{z\in Z}$.
\end{prop}

\section{Generic sequences}\label{SecGenerics}

In the paper \cite{KubFrais} we introduced and studied sequences leading to universal homogeneous structures.
We now adapt the theory to our setting.
As usual, we assume that $\paircats \fK$ is a normed category with the almost amalgamation property.

A sequence $\map{\wek u} \nat {\lcat \fK}$ is \define{\fra}{\fra\ sequence} in $\paircats \fK$ if it satisfies the following two conditions:
\begin{enumerate}
	\item[(U)] For every $\lcat \fK$-object $x$, for every $\eps > 0$, there exist $\ntr$ and a $\rcat \fK$-arrow $\map f x {u_n}$ such that $\mu(f) < \eps$.
\define{}{(U)}
	\item[(A)] Given $\eps > 0$, given a $\lcat \fK$-arrow $\map f {u_n} y$, there exist $m > n$ and a $\rcat \fK$-arrow $\map g y {u_m}$ such that $\mu(g) < \eps$ and $\rho(u_n^m, g \cmp f) < \eps$.
\define{}{(A)}
\end{enumerate}
Recall that we identify a sequence with all of its cofinal subsequences.
It turns out that the definition above is ``correct" because of the almost amalgamation property:

\begin{prop}\label{Perngoweg}
Assume $\paircats \fK$ is a normed category with the almost amalgamation property.
Let $\wek u$ be a sequence in ${\lcat \fK}$. The following conditions are equivalent.
\begin{enumerate}
	\item[$(a)$] $\wek u$ is \fra\ in $\paircats \fK$.
	\item[$(b)$] $\wek u$ has a cofinal subsequence that is \fra\ in $\paircats \fK$.
	\item[$(c)$] Every cofinal subsequence of $\wek u$ is \fra\ in $\paircats \fK$.
\end{enumerate}
\end{prop}

\begin{pf}
Implications (a)$\implies$(c) and (c)$\implies$(b) are obvious.
In fact, the almost amalgamation property is used only for showing that (b)$\implies$(a).

Suppose $M\subs \nat$ is infinite and such that $\wek u \rest M$ is \fra\ in $\paircats \fK$.
Fix $n \in \nat \setminus M$ and fix a $\lcat \fK$-arrow $\map f {u_n} y$.
Fix $\eps > 0$.
Using the almost amalgamation property, we can find $\lcat \fK$-arrows $\map {f'} {u_m} w$ and $\map j y w$ such that $m \in M$, $m > n$ and the diagram
$$\xymatrix{
u_n \ar[rd]_f \ar[rr]^{u_n^m} & & u_m \ar[rd]^{f'} & \\
& y \ar[rr]_j & & w
}$$
is $\eps/2$-commutative.
Since $\wek u \rest M$ is \fra, there is a $\rcat \fK$-arrow $\map g w u_\ell$ with $\ell > m$, $\mu(g) < \eps$, and such that 
the triangle
$$\xymatrix{
u_m \ar[rd]_{f'} \ar[rr]^{u_m^\ell} & & u_\ell \\
& w \ar[ru]_g &
}$$
is $\eps/2$-commutative.
Finally, $\mu(g \cmp j) \loe \mu(g) < \eps$ and $\rho(g \cmp j \cmp f, u_n^\ell) < \eps$, which shows that $\wek u$ satisfies (A).
\end{pf}

The following characterization of a \fra\ sequence will be used later.

\begin{prop}\label{Puhrg}
Let $\paircats \fK$ be a normed category and let $\wek u$ be a sequence in $\lcat\fK$ satisfying (U).
Then $\wek u$ is \fra\ in $\paircats \fK$ if and only if it satisfies the following condition:
\begin{enumerate}
	\item[{\rm(B)}] Given $\eps > 0$, given $\ntr$, given a $\rcat \fK$-arrow $\map f {u_n} y$ with $\mu(f) < +\infty$, there exist $m > n$ and a $\rcat \fK$-arrow $\map g y {u_m}$ such that
$$\mu(g) < \eps \oraz \rho(g \cmp f, u_n^m) < \mu(f) + \eps.$$
\define{}{(B)}
\end{enumerate}
\end{prop}

\begin{pf}
It is obvious that (B) implies (A).
Suppose $\wek u$ is \fra\ and choose $\lcat \fK$-arrows $\map i {u_n} w$ and $\map j y w$ such that $\rho(j \cmp f, i) < \mu(f) + \eps/2$.
Using (A), find $m > n$ and a $\rcat \fK$-arrow $\map k y {u_m}$ such that $\mu(k) < \eps$ and $\rho(k \cmp i, u_n^m) < \eps/2$.
Define $g = k \cmp j$.
Then $\mu(g) \loe \mu(k) + \mu(j) = \mu(k) < \eps$ and
$$\rho(g \cmp f, u_n^m) \loe \rho(k \cmp j \cmp f, k \cmp i) + \rho(k \cmp i, u_n^m) < \mu(f) + \eps/2 + \eps/2 = \mu(f) + \eps,$$
which shows that $\wek u$ satisfies (A).
\end{pf}

\subsection{Separability}

We now turn to the question of existence of a \fra\ sequence.
Let $\paircats \fK$ be a normed category.

A subcategory $\Ef$ of $\lcat \fK$ is \define{dominating}{dominating subcategory} in $\paircats \fK$ if
\begin{enumerate}
	\item[($D_1$)] Every $\lcat \fK$-object has $\rcat \fK$-arrows into $\Ef$-objects of arbitrarily small norm.
\define{}{($D_1$)}
More precisely, given $x \in \ob{\lcat \fK}$, given $\eps > 0$, there exists $\map f x y$ such that $y \in \ob \Ef$ and $\mu(f) < \eps$.
	\item[($D_2$)] Given $\eps > 0$, a $\lcat \fK$-arrow $\map f a y$ such that $a \in \ob \Ef$, there exist a $\rcat \fK$-arrow $\map g y b$ and an $\Ef$-arrow $\map u a b$  such that $\mu(g) < \eps$ and $\rho(g \cmp f, u) < \eps$.
\define{}{($D_2$)}
\end{enumerate}
In some cases, it will be convenient to consider a condition stronger than ($D_2$), namely:
\begin{enumerate}
	\item[($D_2'$)] Given $\eps > 0$, a $\lcat \fK$-arrow $\map f a y$ such that $a \in \ob \Ef$, there exists a $\rcat \fK$-arrow $\map g y b$ such that $\mu(g) < \eps$ and $g \cmp f \in \Ef$.
\define{}{($D'_2$)}
\end{enumerate}
We shall say that $\Ef$ is \define{strongly dominating}{dominating subcategory!-- strongly} in $\paircats \fK$ if it satisfies ($D_1$) and ($D_2'$).
A normed category $\paircats \fK$ is \define{separable}{normed category!-- separable} if there exists a countable $\Ef \subs \lcat \fK$ that is dominating in $\paircats \fK$.

In the next proof, we shall use the simple folklore fact, known as the Rasiowa-Sikorski Lemma: given a directed partially ordered set $\poset$, given a countable family $\ciag D$ of cofinal subsets of $\poset$, there exists an increasing sequence $\ciag p \subs \poset$ such that $p_n \in D_n$ for every $\ntr$. 

\begin{tw}\label{TwExistence}
Let $\paircats \fK$ be a directed normed category with the almost amalgamation property.
The following conditions are equivalent:
\begin{enumerate}
	\item[$(a)$] $\paircats \fK$ is separable.
	\item[$(b)$] $\paircats \fK$ has a \fra\ sequence.
\end{enumerate}
Furthermore, if $\Ef$ is a countable directed dominating subcategory of $\paircats \fK$ with the almost amalgamation property, then there exists a sequence in $\Ef$ that is \fra\ in $\paircats \fK$.
\end{tw}

\begin{pf}
Implication (b)$\implies$(a) is obvious.

Assume $\Ef \subs \lcat \fK$ is countable and dominating in $\paircats \fK$.
Enlarging $\Ef$ if necessary, we may assume that it is directed and has the almost amalgamation property.
We are going to find a \fra\ sequence in $\Ef$, which by ($D_1$) and ($D_2$) must also be \fra\ in $\paircats \fK$.
This will also show the ``furthermore" statement.

Define the following partially ordered set $\poset$:
Elements of $\poset$ are finite sequences in $\Ef$ (i.e. covariant functors from $n < \nat$ into $\Ef$).
The order is end-extension, that is, $\wek x \loe \wek y$ if $\wek y \rest n = \wek x$, where $n = \dom(\wek x)$.

Fix $n,k \in \nat$ and fix an $\Ef$-arrow $\map f a b$.
We define $D_{f,n,k} \subs \poset$ to be the set of all $\wek x \in \poset$ such that $\dom(\wek x) > n$ and the following two conditions are satisfied:
\begin{enumerate}
	\item[(1)] There exists $\ell < \dom(\wek x)$ such that $\Ef(a, x_\ell) \nnempty$.
	\item[(2)] If $a = x_n$ then there exists an $\Ef$-arrow $\map g b {x_m}$ such that $n \loe m < \dom(\wek x)$ and $\rho(g \cmp f, x_n^m) < 1/k$.
\end{enumerate}
Since $\Ef$ is directed and has the almost amalgamation property, it is clear that all sets of the form $D_{f,n,k}$ are cofinal in $\poset$.
It is important that there are only countably many such sets.
Thus, by the Rasiowa-Sikorski Lemma, there exists a sequence
$$\wek u_0 < \wek u_1 < \wek u_2 < \cdots$$
such that for each suitable triple $f,n,k$ there is $r \in \nat$ satisfying $\wek u_r \in D_{f,n,k}$.
It is now rather clear that $\wek u := \bigcup_{\ntr} \wek u_n$ is a \fra\ sequence in $\Ef$ which, by the remarks above, is also a \fra\ sequence in $\paircats \fK$.
\end{pf}

Let us note that the second part of Theorem~\ref{TwExistence} may give some additional information on the structure of the ``approximate" \fra\ limit associated to the \fra\ sequence.
For example, this is the case with the \Gurarii\ space, where there is a countable dominating subcategory whose objects are precisely the $\ell_\infty^n$ spaces, which shows that the \Gurarii\ space is a so-called Lindenstrauss space (see \cite{GarKub1} for more details).

\subsection{Approximate back-and-forth argument}

We now show that a \fra\ sequence is ``almost homogeneous" in the sense described below.

\begin{lm}\label{Labnf}
Let $\paircats \fK$ be a normed category and let $\wek u$, $\wek v$ be \fra\ sequences in $\paircats \fK$.
Furthermore, let $\eps > 0$ and let $\map h {u_0} {v_0}$ be a $\rcat \fK$-arrow with $\mu(h) < \eps$.
Then there exists an approximate isomorphism $\map {H} {\wek u} {\wek v}$ such that $\mu(H) = 0$ and the diagram
$$\xymatrix{
\wek u \ar[rr]^{H} & & \wek v \\
u_0 \ar[u]^{u_0^\infty} \ar[rr]_{h} & & v_0 \ar[u]_{v_0^\infty}
}$$
is $\eps$-commutative.
\end{lm}

\begin{pf}
Fix a decreasing sequence of positive reals $\ciag \eps$ such that
$$
\mu(h) < \eps_0 < \eps \oraz 2 \sum_{n=1}^\infty \eps_n < \eps - \eps_0.
$$
We define inductively sequences of $\rcat \fK$-arrows $\map {f_n} {u_{\phi(n)}} {v_{\psi(n)}}$, $\map {g_n} {v_{\psi(n)}} {u_{\phi(n+1)}}$ such that
\begin{enumerate}
	\item[(1)] $\phi(n) \loe \psi(n) < \phi(n+1)$;
	\item[(2)] $\rho(g_n \cmp f_n, u_{\phi(n)}^{\phi(n+1)}) < \eps_n$;
	\item[(3)] $\rho(f_n \cmp g_{n-1}, v_{\psi(n-1)}^{\psi(n)}) < \eps_n$;
	\item[(4)] $\mu(f_n) < \eps_n$ and $\mu(g_n) < \eps_{n+1}$;
\end{enumerate}
We start by setting $\phi(0) = \psi(0) = 0$ and $f_0 = h$.
We find $g_0$ and $\phi(1)$ by using condition (B) of Proposition~\ref{Puhrg}.

We continue, using condition (B) for both sequences repeatedly.
More precisely, having defined $f_{n-1}$ and $g_{n-1}$, we first use property (B) of the sequence $\wek v$ is \fra, constructing $f_n$ satisfying (3) and with $\mu(f_n) < \eps_n$; next we use the fact that $\wek u$ satisfies (B) in order to find $g_n$ satisfying (2) and with $\mu(g_n) < \eps_{n+1}$.

We now check that $\wek f = \ciag f$ and $\wek g = \ciag g$ are approximate arrows.
Fix $\ntr$ and observe that
\begin{align*}
\rho \Bigl(v_{\psi(n)}^{\psi(n+1)} &\cmp f_n, f_{n+1} \cmp u_{\phi(n)}^{\phi(n+1)} \Bigr) \\
&\loe \rho \Bigl(v_{\psi(n)}^{\psi(n+1)} \cmp f_n, f_{n+1} \cmp g_n \cmp f_n \Bigr) + \rho \Bigl(f_{n+1} \cmp g_n \cmp f_n, f_{n+1} \cmp u_{\phi(n)}^{\phi(n+1)} \Bigr) \\
&\loe \rho \Bigl(v_{\psi(n)}^{\psi(n+1)}, f_{n+1} \cmp g_n \Bigr) + \rho \Bigl(g_n \cmp f_n, u_{\phi(n)}^{\phi(n+1)} \Bigr) < \eps_{n+1} + \eps_n.
\end{align*}
Since the series $\sum_{\ntr}\eps_n$ is convergent, we conclude that $\ciag f$ is an approximate arrow from $\wek u$ to $\wek v$.

By symmetry, we deduce that $\wek g$ is an approximate arrow from $\wek v$ to $\wek u$.
Conditions (2) and (3) tell us that the compositions $\wek f \cmp \wek g$ and $\wek g \cmp \wek f$ are equivalent to the identities, which shows that $H := \wek f$ is an isomorphism.
Condition (4) ensures us that $\mu(H) = 0$.
Finally, recalling that $h = f_0$, we obtain
\begin{align*}
\rho(v_0^\infty \cmp h, H \cmp u_0^\infty) &\loe \sum_{n=0}^\infty \rho \Bigl(v_{\psi(n)}^{\psi(n+1)} \cmp f_n, f_{n+1} \cmp u_{\phi(n)}^{\phi(n+1)} \Bigr) \\
&< \sum_{n=0}^\infty (\eps_n + \eps_{n+1}) = \eps_0 + 2\sum_{n=1}^\infty \eps_n < \eps.
\end{align*}
This completes the proof.
\end{pf}

The lemma above has two interesting corollaries.
Recall that a \emph{0-isomorphism} is an isomorphism $H$ with $\mu(H) = 0$.
In such a case also $\mu(H^{-1}) = 0$.

\begin{tw}[Uniqueness]\label{TWnCorUniq}
A normed category $\paircats \fK$ may have at most one \fra\ sequence, up to an approximate 0-isomorphism.
\end{tw}

\begin{tw}[Almost homogeneity]\label{TwAlmHomgy}
Assume $\paircats \fK$ is a normed category with the almost amalgamation property and with a \fra\ sequence $\wek u$.
Then for every $\lcat \fK$-objects $a,b$, for every approximate 0-arrows $\map i a {\wek u}$, $\map j b {\wek u}$, for every $\rcat \fK$-arrow $\map f a b$, for every $\eps > 0$ such that $\mu(f) < \eps$, there exists an approximate 0-isomorphism $\map {H} {\wek u} {\wek u}$ such that the diagram
$$\xymatrix{
\wek u \ar[r]^{H} & \wek u \\
a \ar[u]^{i} \ar[r]_{f} & b \ar[u]_{j}
}$$
is $\eps$-commutative.
\end{tw}

Note that the existence of a \fra\ sequence automatically implies directedness.

\begin{pf}
Recall that, by definition, $i = \sett{i_n}{n \goe n_0}$, where $\lim_{n \goe n_0} \rho(u_n^\infty \cmp i_n, i) = 0$ and $\lim_{n \goe n_0} \mu(i_n) = 0$.
The same applies to $j$.
Choose $\delta > 0$ such that $\mu(f) < \eps - 6 \delta$.
Choose $k$ big enough so that
\begin{equation}
\rho(u_k^\infty \cmp i_k, i) < \delta \oraz \rho(u_k^\infty \cmp j_k, j) < \delta
\label{Eqertoger}
\end{equation}
holds and $\mu(i_n) < \delta$, $\mu(j_n) < \delta$ whenever $n \goe k$.
Let $f_1 = j_k \cmp f$.
Then $\mu(f_1) \loe \mu(f) + \mu(j_k) < \eps - 5 \delta$.
Using Proposition~\ref{Pamelias}, we find $\lcat \fK$-arrows $\map {f_2}{u_k}w$ and $\map {g_1}{u_k}w$ such that
\begin{equation}
\rho(g_1 \cmp f_1, f_2 \cmp i_k) < \eps - 4 \delta.
\label{Eqneob}
\end{equation}
Using the fact that $\wek u$ is \fra, we find $\ell > k$ and $\map {g_2} w {u_\ell}$ such that
\begin{equation}
\mu(g_2) < \delta \oraz \rho(g_2 \cmp g_1, u_k^\ell) < \delta.
\label{Eqtihggs}
\end{equation}
Define $g = g_2 \cmp f_2$.
Then $\mu(g) < \delta$ and the sequences $\sett{u_n}{n \goe k}$, $\sett{u_n}{n \goe \ell}$ are \fra, therefore by Lemma~\ref{Labnf} there exists an approximate 0-isomorphism $\map H {\wek u}{\wek u}$ satisfying
\begin{equation}
\rho(u_\ell^\infty \cmp g, H \cmp u_k^\infty) < \delta.
\end{equation}
The situation is described in the following diagram
$$\xymatrix{
a \ar[dd]_f \ar[ddr]^{f_1} \ar[r]^{i_k} & u_k \ar[dr]^{f_2} \ar[rrrr] & & & & \cdots \ar[r] & \wek u \ar[dd]^H \\
& & w \ar[dr]^{g_2} & & & & \\
b \ar[r]_{j_k} & u_k \ar[ur]^{g_1} \ar[rr] & & u_\ell \ar[rr] & & \cdots \ar[r] & \wek u
}$$
where the first triangle is commutative, the second one is $\delta$-commutative, and the internal square is $(\eps - 4\delta)$-commutative.
Applying (\ref{Eqertoger}), (\ref{Eqneob}), (\ref{Eqtihggs}), the triangle inequality and inequalities (\ref{Equametric}), we obtain
\begin{align*}
\rho(j \cmp f &, H \cmp i)
\loe \rho(j \cmp f, u_k^\infty \cmp j_k \cmp f) + \rho(u_k^\infty \cmp f_1, u_\ell^\infty \cmp g_2 \cmp g_1 \cmp f_1) \\
&+ \rho(u_\ell^\infty \cmp g_2 \cmp g_1 \cmp f_1, u_\ell^\infty \cmp g_2 \cmp f_2 \cmp i_k) + \rho(u_\ell^\infty \cmp g \cmp i_k, H \cmp u_k^\infty \cmp i_k) \\
&+ \rho(H \cmp u_k^\infty \cmp i_k, H \cmp i) \\
&\loe \rho(j, u_k^\infty \cmp j_k) + \rho(u_k^\ell, g_2 \cmp g_1) + \rho(g_1 \cmp f_1, f_2 \cmp i_k) + \rho(u_\ell^\infty \cmp g, H \cmp u_k^\infty) \\
&+ \rho(u_k^\infty \cmp i_k, i) \\
&< \delta + \delta + (\eps - 4\delta) + \delta + \delta = \eps.
\end{align*}
This completes the proof.
\end{pf}

\subsection{Weak universality}

Below we show that all sequences ``embed" into the \fra\ sequence. In model theory, this is known as \emph{universality}, although when using the language of category theory we should avoid confusion saying that the \fra\ sequence $\wek u$ is \define{weakly terminal}{weakly terminal object}, that is, every other sequence has an approximate 0-arrow into the $\wek u$.

\begin{tw}\label{Thunivweek}
Assume $\paircats \fK$ is a normed category with the almost amalgamation property and with a \fra\ sequence $\wek u$.
Then for every sequence $\wek x$ in $\lcat \fK$ there exists an approximate arrow
$$\map {\wek f} {\wek x} {\wek u}$$ such that $\mu(\wek f) = 0$.
\end{tw}

\begin{pf}
We construct inductively a strictly increasing sequence of natural numbers $\ciag k$ and a sequence of $\rcat \fK$-arrows $\map {f_n} {x_n} {u_{k_n}}$ satisfying for each $\ntr$ the condition
\begin{enumerate}
	\item[($*$)] $\rho(u_{k_n}^{k_{n+1}} \cmp f_n, f_{n+1} \cmp x_n^{n+1}) < 3 \cdot 2^{-n}$ and $\mu(f_n) < 2^{-n}$.
\end{enumerate}
We start by finding $f_0$ using condition (U) of a \fra\ sequence.
Fix $\ntr$ and suppose $f_n$ and $k_n$ have been defined already.

Let $\eps = 2^{-n}$.
Since $\mu(f_n) < 2^{-n}$, there exist $\lcat \fK$-arrows $\map i {x_n} v$, $\map j {u_{k_n}} v$ such that
$$\rho(i, j \cmp f_n) < 2^{-n}.$$
Next, using the almost amalgamation property, we find $\lcat \fK$-arrows $\map k v w$ and $\map \ell {x_{n+1}} w$ such that
$$\rho(k \cmp i, \ell \cmp x_n^{n+1}) < 2^{-n}.$$
Finally, using the fact that $\wek u$ is \fra, find $k_{n+1} > k_n$ and a $\rcat \fK$-arrow $\map g w {u_{k_{n+1}}}$ such that $\mu(g) < 2^{-(n+1)}$ and
$$\rho(g \cmp k \cmp j, u_{k_n}^{k_{n+1}}) < 2^{-n}.$$
The situation is described in the following diagram, where both internal squares and the triangle are $2^{-n}$-commutative:
$$\xymatrix{
u_{k_n} \ar[rr] \ar[dr]_j & &  u_{k_{n+1}} \ar[r] & \cdots \\
 & v \ar[r]^k &  w \ar[u]_g \\
x_n \ar[rr] \ar[uu]^{f_n} \ar[ur]^i & & x_{n+1} \ar[u]_\ell \ar[r] & \cdots
}$$
Define $f_{n+1} := g \cmp \ell$.
Then $\mu(f_{n+1}) \loe \mu(g) + \mu(\ell) = \mu(g) < 2^{-(n+1)}$.
The diagram above shows that condition ($*$) is satisfied.
This completes the inductive construction.

Finally, $\wek f = \ciag f$ is an approximate arrow from $\wek x$ to $\wek u$ satisfying $\mu(\wek f) = \lim_{n\to\infty}\mu(f_n) = 0$.
\end{pf}

\section{Applications}\label{SectAppsy}

In this section we collect selected applications of \fra\ sequences in normed categories.

First of all, summarizing the results of Section~\ref{SecGenerics} we arrive at the following algorithm for finding almost homogeneous structures:
\begin{enumerate}
	\item Choose a directed metric-enriched category $\lcat \fK$ for which there is a chance to have a sequence leading to an almost homogeneous object in some bigger category.
	\item Does $\lcat \fK$ have the almost amalgamation property?
	\item If the answer to Question 2 is negative, STOP.
Otherwise, find a natural bigger category $\rcat \fK$ with the same objects as $\lcat \fK$, such that $\paircats \fK$ becomes a normed category.
	\item Is $\paircats \fK$ separable?
	\item If the answer to Question 4 is negative, STOP.
Otherwise, there is a unique \fra\ sequence $\wek u$ in $\paircats \fK$ which is almost homogeneous with respect to $\lcat \fK$-objects.
	\item Interpret $\ciagilr \fK$ in some natural category, that is, find a canonical ``co-limiting" functor from $\ciagilr \fK$ onto a category $\fC$ containing $\lcat \fK$, in which all countable $\lcat \fK$-sequences have co-limits.
	\item Finally, $\U = \lim \wek u$ is the desired $\fC$-object  that is almost homogeneous with respect to $\lcat \fK$.
The object $\U$ is unique up to a 0-isomorphism and every other $\fC$-object has a 0-arrow into $\U$.
\end{enumerate}
In the next subsections we demonstrate the use of this algorithm for obtaining simpler proofs of the existence and properties of some almost homogeneous structures as well as for finding new ones.

Of course, the prototype example is the category of finite metric spaces with isometric embeddings and non-expansive mappings.
When treated as a normed category, it is indeed directed, separable, and has the strict amalgamation property.
Its \fra\ sequence leads to the well known Urysohn space, which is apparently homogeneous with respect to finite metric spaces.
Isometric uniqueness of the Urysohn space follows easily from the homogeneity.
We skip the details here, referring the readers to the survey article of Melleray~\cite{Mell}.

\subsection{The \Gurarii\ space}\label{SubSectrigeuggurarii}

The \Gurarii\ space is the unique separable Banach space $\G$ satisfying the following condition:
\begin{enumerate}
\item[(G)] Given $\eps > 0$, given finite-dimensional Banach spaces $X \subs Y$, given an isometric embedding $\map f X \G$, there exists an $\eps$-isometric embedding $\map g Y \G$ such that $g \rest X = f$.
\define{}{(G)}
\end{enumerate}
Recall that a linear operator $\map T E F$ is an \emph{$\eps$-isometric embedding} if $(1+\eps)^{-1}\norm x < \norm {Tx} < (1+\eps)\norm x$ holds for every $x \in E$.

The \Gurarii\ space was constructed by \Gurarii~\cite{Gur} in 1966, where it was shown that it is almost homogeneous in the sense that every linear isometry between finite-dimensional subspaces of $\G$ extends to a bijective $\eps$-isometry of $\G$.
Furthermore, an easy back-and-forth argument shows that the \Gurarii\ space is unique up to an $\eps$-isometry for every $\eps > 0$.
The question of uniqueness of $\G$ up to a linear isometry was open for some time, solved by Lusky~\cite{Lusky} in 1976, using rather advanced methods.
The first completely elementary proof of the isometric uniqueness of $\G$ has been found very recently by Solecki and the author~\cite{KubSol}.

It turns out that our framework explains both the existence of $\G$, its isometric uniqueness and almost homogeneity with respect to isometries (already shown in \cite{KubSol}).
Actually, a better understanding of the \Gurarii\ space was one of the main motivations for our study.

Namely, as in Examples~\ref{Exebrigo} and~\ref{Extrgergx}, let $\rcat \fB$ and $\lcat \fB$ be the category of finite-dimensional Banach spaces with linear operators of norm $\loe 1$ and with isometric embeddings, respectively.
For simplicity, we consider real Banach spaces, although the same arguments work for the complex ones.
As mentioned in Example~\ref{Exebrigo}, $\paircats \fB$ is a normed category.

By Proposition~\ref{PrpPuszautBanachs}, $\lcat \fB$ has strict amalgamations.
Obviously, it is directed.
Let $\Ef$ be the subcategory of $\lcat \fB$ whose objects are Banach spaces of the form $\pair {\Err^n}{\anorm}$, where the norm is given by the formula
\begin{equation}
\norm x = \max_{i < k} |\phi_i(x)|
\tag{Q}\label{Eqbieitr}
\end{equation}
in which each $\phi_i$ is a functional satisfying $\img {\phi_i}{\Qyu^n} = \Qyu$.
Call such a space \define{rational}{rational Banach space}.
An arrow of $\Ef$ is a linear isometry $\map f {\Err^n}{\Err^m}$ satisfying $\img f{\Qyu^n} \subs \Qyu^m$.
It is clear that $\Ef$ is countable.

\begin{lm}
$\Ef$ is dominating in $\paircats \fB$.
\end{lm}

\begin{pf}
It is rather clear that $\Ef$ satisfies ($D_1$).
In order to see ($D_2$), fix an isometric embedding $\map f X Y$, where $X$ is a rational Banach space.
We may assume that $X = \Err^n$ and $Y = X \oplus \Err^k$ and $f(x) = \pair x0$ for $x \in X$.
Given $\eps > 0$, there exist functionals $\phi_0, \dots, \phi_m$ on $Y$ such that $\norm y_Y$ is $\eps/2$-close to
$$\norm y' = \max_{i < m}|\phi_i(y)|.$$
We may assume that $\norm {\phi_i} \loe 1$ for each $i < m$ and that some of the $\phi_i$s are extensions of the rational functionals defining the norm on $X$.
We can now ``correct" each non-rational $\phi_i$ so that it becomes rational and the respective norm is $\eps/2$-close to $\anorm'$.
Finally, the same map $f$ becomes an isometric embedding $\eps$-close to the original one, when $Y$ is endowed with the new norm.
This shows ($D_2$).
\end{pf}

Thus, $\paircats \fB$ is separable, therefore it has a \fra\ sequence.
Example~\ref{Extrgergx} shows that $\ciagilr \fB$ has a canonical co-limiting functor onto the category of separable Banach spaces.
We still need to translate condition (G).
By a chain of Banach spaces we mean a chain $\ciag X$, where each $X_n$ is a Banach space and the norm of $X_{n+1}$ extends the norm of $X_n$ for every $\ntr$.

\begin{lm}\label{Leenrfi}
Let $\ciag G$ be a chain of finite-dimensional Banach spaces with $G_\infty = \ovr{\bigcup_{\ntr}G_n}$.
If $G_\infty$ satisfies (G) then $\ciag G$ is a \fra\ sequence in $\paircats \fB$.
\end{lm}

\begin{pf}
We check that $\wek G = \ciag G$ satisfies (A).
Fix $\eps > 0$ and an isometric embedding $\map f {G_n} Y$, where $Y$ is a finite-dimensional Banach space.
Using condition (G) for the map $\map {f^{-1}}{\img f{G_n}}{G_\infty}$, we find an $\eps$-isometric embedding $\map h Y {G_\infty}$ such that $h(f(x)) = x$ for every $x \in G_n$.
Using the fact that $Y$ is finite-dimensional, we can find $m > n$ and an $\eps$-isometric embedding $\map {h_1} Y {G_m}$ that is $\eps$-close to $h$.
Finally, define $g := (1+\eps)^{-1} h_1$.
Then $\map g Y {G_m}$ is a $\rcat \fB$-arrow, because $\norm g \loe 1$.
Clearly, $g$ is $\eps$-close to $h_1$.
Finally, $g$ is $2\eps$-close to $h$, therefore $\norm{g(f(x)) - x} = \norm{g(f(x)) - h(f(x))} \loe 2\eps \norm x$ for $x \in G_n$.
This shows (A).

Since $\rcat \fB$ has the initial object $\sn 0$, condition (U) follows from (A).
\end{pf}

It turns out that the converse to the lemma above is also true, because of the uniqueness of the \fra\ sequence, up to a linear isometry of its co-limit:

\begin{lm}
Let $\wek x$, $\wek y$ be sequences in $\lcat \fB$ and let $\map {\wek f}{\wek x}{\wek y}$ be an approximate 0-arrow in $\ciagilr \fB$.
Let $X$ and $Y$ be the co-limits of $\wek x$, $\wek y$ in the category of Banach spaces.
Then $\lim \wek f$ is an isometric embedding of $X$ into $Y$.
\end{lm}

\begin{pf}
We may assume that $\wek x = \ciag X$, $\wek y = \ciag Y$ are chains of (finite-dimensional) Banach spaces and that $\map {f_n}{X_n}{Y_n}$ for each $\ntr$.
Since $\mu(f_n) \to 0$, by Example~\ref{Exebrigo}, we know that given $\eps>0$ there is $n_0$ such that
$$(1 - \eps) \norm x \loe \norm {f_n(x)} \loe \norm x$$
holds for every $n \goe n_0$ and for every $x \in X_n$.
Since $\ciag f$ is an approximate arrow, for each $x \in X_k$, the limit
$$f(x) = \lim_{n > k} f_n(x)$$
exists, because $\sett{f_n(x)}{n > k}$ is a Cauchy sequence.
Thus $f(x)$ is defined on $\bigcup_{\ntr}X_n$, therefore it has a unique extension to a continuous linear operator from $X$ to $Y$ which is an $\eps$-isometric embedding for every $\eps > 0$. Thus, the extension of $f$ is an isometric embedding.
\end{pf}

Thus, a \fra\ sequence in $\paircats \fB$ yields an isometrically unique separable Banach space $\G$, which must be the \Gurarii\ space by Lemma~\ref{Leenrfi}.
As we have mentioned before, an elementary proof of its isometric uniqueness \cite{KubSol} (using Example~\ref{Exebrigo}) was one of the main inspirations for studying \fra\ sequences in the context of metric-enriched categories.

\separator

We shall now describe a natural category leading to a universal projection.
This will bring an improvement of a result due to Wojtaszczyk~\cite{Wojt} and Lusky~\cite{Lusky2} on complemented subspaces of the \Gurarii\ space.
In particular, we shall show that there exists a projection $P$ of the \Gurarii\ space $\G$ whose range and kernel are linearly isometric to $\G$ and $P$ is almost homogeneous with respect to finite-dimensional subspaces and ``contains" all norm one operators between separable Banach spaces.
The construction is inspired by a recent work of Pech \& Pech~\cite{PechPech} involving comma categories.
Recall that a Banach space $X$ is \define{complemented}{complemented Banach space} in a space $Y$ if $X \subs Y$ and there is a bounded linear operator $\map P Y X$ such that $P \rest X = \id X$; such an operator is called a \define{projection}{projection}.
A space $X$ is \emph{$1$-complemented} in $Y$ if there is a projection of norm 1 witnessing that $X$ is complemented in $Y$.

From now on, we fix a separable Banach space $\bS$.
We shall define a category $\rcat {\fK(\bS)}$ as follows.
The objects will be linear operators of the form $\map T E \bS$ satisfying $\norm T \loe 1$ and such that $E$ is a finite-dimensional space.
An arrow from $\map T E \bS$ to $\map {T'}{E'}\bS$ will be a linear operator $\map f E {E'}$ such that $\norm f \loe 1$.
We shall say that $f$ is a $\lcat {\fK(\bS)}$-arrow if it is an isometric embedding and $T' \cmp f = T$.
There is an obvious metric on $\rcat {\fK(\bS)}$, namely $\rho(f,g) =\norm{f-g}$, whenever $f$ and $g$ are $\rcat {\fK(\bS)}$-arrows from $T$ to $T'$.
Note that $\lcat {\fK(\bS)}$ is a comma category based on $\bS$, restricted to finite-dimensional spaces.

It is important to ``decode" the norm given by the pair $\paircats {\fK(\bS)}$.
This is given below.

\begin{lm}\label{Lmytvhfrut}
Let $E$, $F$ be finite-dimensional Banach spaces, let $\map T E Y$, $\map R F Y$ be linear operators of norm $\loe1$ and let $\map f E F$ be an injective linear operator such that $\norm f \loe 1$, $\norm {f^{-1}} \goe 1-\eps$ and $\norm {R \cmp f - T} \loe \eps$, where $\eps > 0$ is fixed.
Then there exist isometric embeddings $\map i E Z$, $\map j F Z$ and a linear operator $\map U Z Y$, where $Z$ is finite-dimensional, $U \cmp i = T$, $U \cmp j = R$, and $$\norm { j \cmp f - i } \loe \eps.$$

In other words, $\mu(f) \loe \eps$, with respect to $\paircats {\fK(\bS)}$.
\end{lm}

\begin{pf}
Looking at Example~\ref{Exebrigo} above, it is not hard to show the following property:
\begin{enumerate}
	\item[(*)] Given linear operators $\map p E V$, $\map q F V$ of norm $\loe1$ and such that $\norm {p - q \cmp f} \loe \eps$, the unique linear operator $S$ from $Z = E \oplus F$ with the norm defined in Example~\ref{Exebrigo} into $V$, satisfying $S \cmp i = p$ and $S \cmp j = q$, has norm $\loe1$. 
\end{enumerate}
The details are given in~\cite{JoasiaG}.
It is obvious that property (*) implies the lemma.
\end{pf}

\begin{lm}\label{Lmbeigheui}
The pair $\pair {\lcat {\fK(\bS)}}{\rcat {\fK(\bS)}}$ is a separable directed normed category such that $\lcat {\fK(\bS)}$ has the strict amalgamation property.
\end{lm}

\begin{pf}
It is clear that $\lcat {\fK(\bS)}$ is directed and has the strict amalgamation property, the latter follows from the fact that the category of Banach spaces admits pushouts (see Proposition~\ref{PrpPuszautBanachs} above).
It remains to show that it is separable.

Fix a countable dense $\Qyu$-linear subspace $\bS_0$ of $\bS$ and define $\Ef$ to be the family of all $\lcat {\fK(\bS)}$-arrows $f$ from $\map T E \bS$ to $\map {T'}{E'}\bS$, where $E$ and $E'$ are rational Banach spaces (i.e. $E = \Err^n$, $F = \Err^m$ and the norms are induced by finitely many rational functionals, see formula (\ref{Eqbieitr}) above), the operators $T$, $T'$ map rational vectors into $\bS_0$, and $f$ maps rational vectors to rational vectors (a vector in $\Err^k$ is \define{rational}{rational vector} if its coordinates are rational).
Obviously, $\Ef$ is countable and it is easy to check (using Lemma~\ref{Lmytvhfrut}) that it is dominating in $\pair {\lcat {\fK(\bS)}}{\rcat {\fK(\bS)}}$.
\end{pf}

One has to stress out the importance of Lemma~\ref{Lmytvhfrut}.
Without it, we would only know that $\paircats {\fK(\bS)}$ is a normed category, without having any description of its norm.
In the extreme case, a normed category $\paircats \fK$ can have the property that $\mu(f) = +\infty$ whenever $f \in \rcat \fK \setminus \lcat \fK$.

By Lemma~\ref{Lmbeigheui}, the normed category $\paircats {\fK(\bS)}$ has a \fra\ sequence $\sett{\map {U_n}{u_n} \bS}{\ntr}$.
Without loss of generality, we may assume that $u_n \subs u_{n+1}$ for each $\ntr$.
Denote by $\G(\bS)$ the completion of the chain $\ciag u$ and let $\map {U_\bS}{\G(\bS)} \bS$ be the unique linear operator satisfying $U_\bS \rest u_n = U_n$ for every $\ntr$.

The next result, showing the properties of $U_\bS$, is a straightforward application of Theorems~\ref{TwAlmHomgy} and \ref{Thunivweek}.

\begin{tw}
The operator $\map {U_\bS}{\G(\bS)} \bS$ has the following properties.
\begin{enumerate}
	\item[(1)] For every linear operator $\map T X \bS$ of norm $\loe1$, with $X$ separable, there exists an isometric embedding $\map j X {\G(\bS)}$ such that $U_\bS \cmp j = T$.
	\item[(2)] Given $\eps>0$, given finite-dimensional spaces $E, E' \subs \G(\bS)$, given an $\eps$-isometric embedding such that $\norm{U_\bS \cmp f - U_\bS \rest E} < \eps$, there exists a bijective linear isometry $\map h {\G(\bS)}{\G(\bS)}$ such that $U_\bS \cmp h = U_\bS$ and $\norm{h \rest E - f} < \eps$.
	\item[(3)] $U_\bS$ is right-invertible and its kernel is linearly isometric to the \Gurarii\ space.
	\item[(4)] Conditions (1) and (2) determine $U_\bS$ uniquely, up to a linear isometry.
\end{enumerate}
\end{tw}

\begin{pf}
Property (1) follows from Theorem~\ref{Thunivweek}, knowing that every operator 
is the co-limit of a sequence in $\lcat{\fK(\bS)}$.
Property (2) follows from Theorem~\ref{TwAlmHomgy}, having in mind Lemma~\ref{Lmytvhfrut}.
The fact that $U_\bS$ is left-invertible follows from (1) applied to the identity $\map {\id\bS} \bS \bS$.
Let $G = \ker U_\bS$.
Notice that every operator $\map h {\G(\bS)}{\G(\bS)}$ satisfying $U_\bS \cmp h = U_\bS$ preserves $G$.
Thus, property (2) applied to subspaces of $G$ shows that $G$ is almost homogeneous with respect to its finite-dimensional subspaces.
Property (1) applied to zero operators $\map {0_X} X \bS$ shows that $G$ contains isometric copies of all separable spaces, therefore it is linearly isometric to $\G$.
Finally, uniqueness of $U_\bS$ follows from the fact that any of its decomposition into a chain of operators on finite-dimensional spaces leads to a \fra\ sequence in $\paircats {\fK(\bS)}$.
\end{pf}

\begin{uwgi}
The work \cite{GarKub0} contains a construction of a universal linear operator on the \Gurarii\ space.
Again, it is possible to describe it in the language of normed categories.
Namely, the objects of the category $\rcat \fK$ are linear operators $\map T {X_0} {X_1}$, where $X_0$, $X_1$ are finite-dimensional Banach spaces and $\norm T \loe 1$. An arrow from $T$ to $S$ is a pair $\pair {f_0}{f_1}$ of linear operators of norm $\loe 1$ satisfying $S \cmp f_0 = f_1 \cmp T$.
Obviously, $\lcat \fK$ should be the subcategory of all pairs of isometric embeddings.
The main lemma in~\cite{GarKub0} says that $\paircats \fK$ has the strict amalgamation property.
The remaining issues are easily solved and the \fra\ sequence in $\paircats \fK$ leads to the universal (almost homogeneous) linear operator whose domain and range turn out to be isometric to the \Gurarii\ space.
The details can be found in~\cite{GarKub0}, actually without referring to normed categories.
\end{uwgi}

\subsection{The Cantor set}

We are going to revisit some folklore facts about the Cantor set, the simplest reversed \fra\ limit in the sense of~\cite{IrwSol}.
Some of our ideas are already contained in \cite{FelixC}.

Namely, let $\rcat \fK$ be the opposite category of all non-expansive maps between non\-empty compact metric spaces and let $\lcat \fK \subs \rcat \fK$ be the opposite category of all quotient maps.
Formally, a $\rcat \fK$-arrow $f$ from $X$ to $Y$ is a non-expansive map $\map f Y X$.
In particular, a \emph{sequence} in $\lcat \fK$ is an inverse sequence of compact metric spaces in which all bonding maps are non-expansive.

It is clear that $\rcat \fK$ is metric-enriched, by setting
$$\rho(f,g) = \max_{t\in K} d(f(t),g(t)),$$
where $\map {f,g}K L$ and $d$ is the metric of $L$.

The following simple fact will be needed later.

\begin{lm}\label{Lmergerol}
Let $\pair X d$ be a metric space and let $\sett{\map {f_i}X {Y_i}}{i<n}$ be a finite family of continuous maps into metric spaces $\pair {Y_i}{d_i}$.
Then the formula
$$\rho(s,t) = d(s,t) + \max_{i<n}d_i(f_i(s), f_i(t))$$
defines a compatible metric on $X$ such that all maps $f_i$ become non-expansive.
\end{lm}

It is well known that given two quotient maps $\map f X Z$, $\map g Y Z$, there exist quotient maps $\map {f'} W X$, $\map {g'} W Y$ such that $f \cmp f' = g \cmp g'$.
In fact, one may define
$$W = \setof{\pair st \in X \times Y}{f(s) = g(t)}$$
endowed with the product topology, and $f'$, $g'$ are the projections.
Actually, this is the standard construction of the pushout in the category of topological spaces.
Now, if $Z,X,Y$ are metric spaces and $f$, $g$ are non-expansive, then Lemma~\ref{Lmergerol} gives a metric on $W$ for which $f'$, $g'$ become non-expansive.
This shows that the category $\lcat \fK$ has the strict amalgamation property.
Since the singleton is initial, $\lcat \fK$ is directed.

We now describe the norm $\mu$ on $\paircats \fK$.

\begin{prop}\label{Prpxbewhu}
Let $\map f X Y$ be a non-expansive map of compact metric spaces.
Then
$$\mu(f) = \inf \setof{ \eps > 0 }{ \img f X \text{ is $\eps$-dense in }Y }.$$
\end{prop}

\begin{pf}
Suppose $\mu(f) < \eps$ and choose 1-Lipschitz quotient maps $\map p Z Y$, $\map q Z X$ such that $\rho(p, f \cmp q) < \eps$.
Fix $y \in Y$ and choose $z \in Z$ such that $y = p(z)$.
Then $\eps > \rho(p, f \cmp q) \goe d_Y(p(z), f(q(z)))$.
Since $f(q(z)) \in \img {f \cmp q} Z = \img f X$, we conclude that $\dist(y, \img f X) < \eps$.

Now suppose that $\img f X$ is $\eps$-dense in $Y$ and define
$$Z = \setof{\pair x y \in X \times Y}{d_Y(y, f(x)) \loe \eps}.$$
Then $Z$ is a closed subspace of $X \times Y$.
Endow it with the maximum metric.
Let $\map p Z Y$, $\map q Z X$ be the projections restricted to $Z$.
It is obvious that $q$ is onto.
Given $y \in Y$, there is $x \in X$ such that $d_Y(y, f(x)) < \eps$, therefore $\pair x y \in Z$ and $y = p(x,y)$.
This shows that $p$ is onto.
Finally, $d_Y(f(q(x,y)), p(x,y)) = d_Y(f(x), y) \loe \eps$, therefore $\rho(f \cmp q, p) \loe \eps$, showing that $\mu(f) \loe \eps$.
\end{pf}

Let $\Ef$ be the subcategory of $\lcat \fK$ consisting of 1-Lipschitz quotient maps between finite rational metric spaces.
Recall that a metric space $\pair X d$ is \define{rational}{rational metric space} if $\img d{X \times X} \subs \Qyu$.
It is clear that $\Ef$ is countable, directed and has the strict amalgamation property (the metric given by Lemma~\ref{Lmergerol} is rational if all the involved metrics are rational).
We are going to show that $\Ef$ is dominating in $\paircats \fK$.
Below is the crucial claim.

\begin{lm}\label{Lemkrcjata}
Let $X$ be a finite metric space and let $\eps > 0$.
Then there exist a finite rational metric space $Y$ and a $1$-Lipschitz bijection $\map h Y X$ such that $\Lip {h^{-1}} < 1 + \eps$.
\end{lm}

\begin{pf}
Fix $\eta > 0$ and define $r = \min\setof{d(s,t)}{s,t \in X, \; s \ne t}$.
Let $\delta = \eta \cdot r / 2$.
We may assume that $X$ is a subspace of the Urysohn space $\U$.
Recall that the rational Urysohn space $\U_\Qyu$ is dense in $\U$, therefore for each $s \in X$ we can find $y_s \in \U_\Qyu$ such that $d(s, y_s) < \delta$.
Let $Y = \setof{y_s}{s\in X}$ and let $\map h Y X$ be the obvious bijection, i.e., $h(y_s) = s$ for $s \in X$.
The space $Y$ is rational, although $h$ may not be 1-Lipschitz.
We shall later enlarge the metric of $Y$ so that $h$ will become 1-Lipschitz.
First, note that given $s,t \in X$ we have
$$d(h(y_s), h(y_t)) = d(s,t) < d(y_s, y_t) + 2\delta = d(y_s, y_t) + 2 \eta \cdot r \loe (1 + \eta) d(y_s, y_t).$$
Similarly, $d(h(y_s), h(y_t)) > (1 - \eta) d(y_s, y_t)$.
It follows that $\Lip {h} < 1 + \eta$ and $\Lip {h^{-1}} < (1 - \eta)^{-1}$.

Now, suppose that $\eta$ is rational and satisfies $(1+\eta)(1-\eta)^{-1} \loe 1 + \eps$.
Consider $Y$ with the metric $d' = (1+\eta)d$.
This is still a rational metric space and $h$ is 1-Lipschitz with respect to $d'$.
Finally, $\Lip {h^{-1}} < (1+\eta)(1-\eta)^{-1} \loe 1 + \eps$ with respect to the metric $d'$.
\end{pf}

\begin{prop}
$\Ef$ is strongly dominating in $\paircats \fK$.
\end{prop}

\begin{pf}
Condition ($D_1$) follows directly from Lemma~\ref{Lemkrcjata}.
Fix $\eps > 0$ and a finite rational metric space $X$ and fix a non-expansive quotient map $\map f K X$, where $K$ is a compact metric space.
Choose a finite $\eps$-dense subset of $K$ such that $\img f S = X$.
Using Lemma~\ref{Lemkrcjata}, find a rational metric space $Y$ and a non-expansive bijection $\map h Y S$ such that $\Lip{h^{-1}} < 1 + \eps$.
Now consider $h$ as a map from $Y$ to $K$.
By Proposition~\ref{Prpxbewhu}, $\mu(h) \loe \eps$.
Finally, $f \cmp h$ is a quotient map of finite rational metric spaces, therefore $f \cmp h \in \Ef$.
\end{pf}

\begin{wn}
$\paircats \fK$ is a separable directed normed category with the strict amalgamation property.
\end{wn}

By Theorem~\ref{TwExistence} there exists a sequence $\wek u$ in $\Ef$ that is \fra\ in $\paircats \fK$.

We shall now get rid of the metrics, moving to the category of compact metric spaces.
More precisely, define the ``co-limiting" functor $\map \lim {\ciagilr \fK}\Komp$ in the obvious way: $\lim \wek x$ should be the inverse limit of the sequence $\wek x$ in the category of topological spaces.
In particular, the functor $\lim$ forgets the metric structures of the sequence.
Notice that, given any inverse sequence $\wek X$ of nonempty compact metrizable spaces with quotient maps, by Lemma~\ref{Lmergerol} and induction, there exist compatible metrics on each $X_n$ such that all bonding maps become 1-Lipschitz.
In fact, Lemma~\ref{Lmergerol} implies more: Given inverse sequences $\wek K$, $\wek L$ of nonempty compact metrizable spaces (with quotient bonding maps), given a natural transformation $\map {\wek f}{\wek K}{\wek L}$, there exist compatible metrics such that all mappings in the diagram
$$\xymatrix{
K_0 \ar[d]_{f_0} & K_1 \ar[l] \ar[d]^{f_1} & \cdots \ar[l] & K_n \ar[d]_{f_n} \ar[l] & K_{n+1} \ar[l] \ar[d]^{f_{n+1}} & \cdots \ar[l] \\
L_0 & L_1 \ar[l] & \cdots \ar[l] & L_n \ar[l] & L_{n+1} \ar[l] & \cdots \ar[l]
}$$
become 1-Lipschitz.

What is more important, every approximate arrow of sequences ``converges" to a continuous map of their inverse limits.
Moreover, an approximate 0-arrow ``converges" to a quotient map of the limits.
This can be checked easily (see also~\cite{Miod1}).

Let $C$ be the inverse limit of $\wek u$ in the category of topological spaces.
As one can expect, $C$ is the Cantor set.
Indeed, $C$ is 0-dimensional, being the inverse limit of finite sets.
Furthermore, $C$ has the following property:
\begin{enumerate}
	\item[(C)] Given a quotient map of nonempty finite sets $\map f T S$, given a quotient map $\map p C S$, there exists a quotient map $\map q C T$ such that $f \cmp q = p$.
\define{}{(C)}
\end{enumerate}
Indeed, assuming $p$ is given, we find $\ntr$ such that $p = p' \cmp u_n^\infty$, where $\map {u_n^\infty} C {u_n}$ is the canonical projection.
Since $u_n$ is a finite metric space, we find $r > 0$ such that $d(s,t) > r$ whenever $s,t \in u_n$ are distinct.
Define a metric $d$ on $S$ by setting $d(x,y) = r$ iff $x \ne y$.
Define the same metric on $T$.
By this way, both $p'$ and $f$ become 1-Lipschitz and we find $q$ by using condition (A) of the \fra\ sequence.

Finally, it is well known and easy to check that a compact space satisfying (C) is dense-it-itself, therefore $C$ is homeomorphic to the Cantor set.

As an application of Theorem~\ref{Thunivweek}, we obtain the folklore fact that every compact metric space is a quotient of the Cantor set.
Below is the translation of almost homogeneity:

\begin{tw}\label{Thhtrght}
Let $K$ be a compact metric space and let $\map p \Cantor K$, $\map q \Cantor K$ be quotient maps.
Then for every $\eps > 0$ there exists a homeomorphism $\map h \Cantor \Cantor$ such that
$$\rho(q \cmp h, p) < \eps.$$
\end{tw}

It is natural to ask whether $\eps$ is needed in the statement above.
The answer is negative, as the following example shows.
Namely, let $\map p \Cantor \unii$ be a quotient map such that $p^{-1}(0)$ is a singleton and $p^{-1}(1)$ contains more than one point.
Let $q = \phi \cmp p$, where $\phi(t) = 1 - t$.
Then there is no homeomorphism $h$ satisfying $q \cmp h = p$.

Finally, let us note that condition (C) characterizes the Cantor set among 0-dimen\-sio\-nal compact metric spaces only.
Indeed, take $C = \Cantor \times \unii$ and observe that for every quotient map $\map f C S$ with $S$ finite, there exists a unique map $\map g \Cantor S$ such that $f = g \cmp p$, where $\map p C \Cantor$ is the canonical projection.
It follows that $C$ satisfies condition (C), yet $C \not\iso \Cantor$.
It turns out however that the condition in Theorem~\ref{Thhtrght} characterizes the Cantor set:

\begin{tw}
Assume $C$ is a compact metrizable space that maps onto all polyhedra and satisfies the assertion of Theorem~\ref{Thhtrght} for every polyhedron $K$.
Then $C$ is homeomorphic to the Cantor set.
\end{tw}

\begin{pf}
By Freudenthal's theorem~\cite{Freud}, $C$ is the inverse limit of a sequence of polyhedra $\ciag \Delta$ with quotient maps $\map {\delta_n^m}{\Delta_m}{\Delta_n}$ ($n < m < \nat$).
By Theorem~\ref{TWnCorUniq}, it suffices to show that $\ciag \Delta$ is a \fra\ sequence in $\paircats \fK$.

Fix $\eps > 0$, $\ntr$ and fix a quotient map $\map f K {\Delta_n}$.
By assumption, there exists a quotient map $\map g C K$.
Using the condition of Theorem~\ref{Thhtrght} for the maps $\delta_n^\infty$ and $f \cmp q$, we find a homeomorphism $\map h C C$ such that
$$\rho(f \cmp q \cmp h, \delta_n^\infty) < \eps/2.$$
The approximation lemma of Eilenberg \& Steenrod~\cite{EilSte} (see also \cite{Miod1}) says that there exist $m > n$ and a map $\map g {\Delta_m} K$ such that
$$\rho(g \cmp \delta_m^\infty, q \cmp h) < \eps/2.$$
In particular, $\img g{\Delta_m}$ is $\eps/2$-dense in $K$ and hence $\mu(g) < \eps$.
Using the inequalities above, we get
$$\rho(\delta_n^m, f \cmp g) = \rho(\delta_n^\infty, f \cmp g \cmp \delta_m^\infty) \loe \rho(\delta_n^\infty, f \cmp q \cmp h) + \rho(f \cmp q \cmp h, f \cmp g \cmp \delta_m^\infty) < \eps.$$
This shows that $\ciag \Delta$ is \fra\ and completes the proof.
\end{pf}

\subsection{The pseudo-arc}\label{SubSectpseudoarc}

We now describe the universal chainable continuum, known under the name \define{pseudo-arc}{pseudo-arc}, as the limit of a \fra\ sequence in a suitable metric category.
Actually, we shall work in a normed category of the form $\pair \fK \fK$, so in particular $\mu = 0$ and only the metric $\rho$ is relevant.

Recall that a \define{continuum}{continuum} is a compact connected metrizable space.
A continuum is \define{chainable}{continuum!-- chainable} (also called \define{snake-like}{continuum!-- snake-like}) if it is homeomorphic to the limit of an inverse sequence of quotient maps of the unit interval.
In particular, a chainable continuum maps onto the unit interval and therefore cannot be degenerate.

It turns out that we can restrict attention to piece-wise linear maps:

\begin{prop}\label{Ppicewiselins}
Every inverse sequence of quotient maps of the unit interval is equivalent to an inverse sequence of piece-wise linear quotient maps of the unit interval.
\end{prop}

\begin{pf}
Let $\wek f = \sett{f_n^m}{n<m<\nat}$ be an inverse sequence of quotient maps of $\unii$.
We construct inductively piece-wise linear quotient maps $\map {g_n^{n+1}}\unii \unii$ such that
$$\rho(f_n^{n+1}, g_n^{n+1}) < 2^{-n}/k_n$$
where
$$k_n = \Lip {g_0^1} \cdot \Lip {g_1^2} \cdot \dots \cdot \Lip {g_{n-1}^n}.$$
Note that every piece-wise linear map is Lipschitz, therefore $k_n$ is well defined.
We set $g_i^j = g_i^{i+1} \cmp \dots \cmp g_{j-1}^j$.
Finally, $\wek g = \sett{g_n^m}{n<m<\omega}$ is an inverse sequence, easily seen to be equivalent to $\wek f$.
\end{pf}

The following fact can be proved easily, using the linear structure of the unit interval:

\begin{lm}\label{Lewfkt}
Let $\wek q = \sett{q_n^m}{n<m<\nat}$ be an inverse sequence of quotient maps of the unit interval with $K = \liminv \wek q$ in the category of compact spaces.
Denote by $\map{q_n} K \unii$ the canonical projection onto the $n$th element of the sequence.
Given a quotient map $\map f K \unii$, for every $\eps > 0$ there exist $\ntr$ and a piece-wise linear quotient map $\map g \unii \unii$ such that
$$\rho(g \cmp q_n, f) < \eps,$$
where $\rho$ is the maximum metric on the space of continuous functions.
\end{lm}

Let $\lcat \fI$ be the opposite category of non-expansive piece-wise linear quotient mappings of the form
$$\map f {\pair \unii {d_0}} {\pair \unii {d_1}},$$
where $\unii$ is the unit interval and $d_i(s,t) = m_i|s-t|$ for some integer constant $m_i > 0$.
Let $\rcat \fI = \lcat \fI$, that is, we shall really work in a single category $\fI$ and all arrows, including approximate arrows of sequences are 0-arrows.
In other words, the objects of $\fI$ are pairs of the form $\pair \unii d$, where $d$ is the usual metric multiplied by a positive integer constant, needed only for making the maps 1-Lipschitz.
In particular, there are only countably many objects.

We could have defined $\fI$ equivalently by saying that its objects are intervals of the form $[0,n]$ with $\Ntr$ endowed with the usual metric, and arrows are 1-Lipschitz quotient maps.

It is clear that $\fI$ is a metric-enriched category, with the same metric $\rho$ as in the case of all nonempty compact metric spaces.
It is also clear that $\fI$ is directed.
The almost amalgamation property follows from the following well known result, sometimes called the \define{uniformization principle}{uniformization principle}:

\begin{tw}[Mountain Climbing Theorem]
Let $\map {f,g}\unii \unii$ be quotient maps that are piece-wise monotone and satisfy $f(i)=i=g(i)$ for $i=0,1$.
Then there exist quotient maps $\map {f',g'}\unii \unii$ such that $f \cmp f' = g \cmp g'$.
\end{tw}

The result above goes back to Homma~\cite{Homma}; the formulation (actually involving finitely many piece-wise monotone quotient maps) is due to to Sikorski \& Zaran\-kie\-wicz~\cite{SikZar}. Another version, involving two functions that are constant on no open subintervals of $\unii$ is due to Huneke~\cite{Huneke}.

The Mountain Climbing Theorem is usually stated for functions $f,g$ satisfying $f(0) = 0 = g(0)$ and $f(1) = 1 = g(1)$, whose graphs can therefore be interpreted as two slopes of the same mountain.
The functions $f',g'$ can be interpreted as the existence of two ``mountain climbings" on these slopes with the property that at each moment of time the travelers have the same altitude (sometimes one of the travelers has to go backwards).
This justifies the name of the theorem.
It is rather clear that given a quotient map $\map f \unii \unii$, there exists a piece-wise linear quotient map $\map {f_1} \unii \unii$ such that $f_1(f(i)) = i$ for $i = 0,1$.
As a corollary, we get:

\begin{prop}
The category $\fI$ has the strict amalgamation property.
\end{prop}

Let us note that the Mountain Climbing Theorem fails for arbitrary quotient maps of the unit interval; an example was first found by Minagawa (quoted in Homma~\cite{Homma}) and independently by Sikorski \& Zarankiewicz~\cite{SikZar}.

As one can easily guess, $\fI$ is separable.
A natural countable dominating subcategory is described below.

We say that a quotient map $\map f \unii \unii$ is \define{rational}{rational quotient} if $f(0) = 0$, $f(1) = 1$, and there is a decomposition $0 = t_0 < t_1 < \dots < t_n = 1$ such that $\sett{t_i}{i < n} \subs \Qyu$, $f \rest [t_i,t_{i+1}]$ is linear for each $i < n$ and $\sett{f(t_i)}{i < n} \subs \Qyu$.
Finally, define $\Ef \subs \fI$ to be the category of all rational quotient maps.
The following fact is rather obvious.

\begin{prop}
The category $\Ef$ is countable and dominating in $\fI$.
\end{prop}

It is clear that the category $\fI$ has a canonical ``limiting" functor, which assigns the inverse limit to a sequence.
From now on, let $\wek u$ be a \fra\ sequence in $\ciagi \fI$ and let
$$P = \liminv \wek u$$
in the category of compact spaces.
It turns out that $P$ is the pseudo-arc.
We explain the details below.

First of all, recall that formally the \define{pseudo-arc}{pseudo-arc} is defined to be a hereditarily indecomposable chainable continuum, which by a result of Bing~\cite{Bing} is known to be unique.
Recall that $K$ is \define{indecomposable}{continuum!-- indecomposable} if it cannot be written as $A \cup B$ where $A,B$ are proper subcontinua.
A continuum $K$ is \define{hereditarily indecomposable}{continuum!--hereditarily indecomposable} if every subcontinuum of $K$ is indecomposable.

\begin{lm}\label{Legrdin}
Let $\wek v = \sett{v_n^m}{n<m<\nat}$ be an inverse sequence of quotient maps of the unit interval.
Then $\liminv \wek v$ is homeomorphic to $P$ if and only if $\wek v$ satisfies the following condition:
\begin{enumerate}
	\item[{\rm(\P)}] Given $\ntr$, $\eps > 0$, given a quotient map $\map f \unii \unii$, there exist $m > n$ and a quotient map $\map g \unii \unii$ such that $\rho(f \cmp g, v_n^m) < \eps$.
\define{}{(\P)}
\end{enumerate}
\end{lm}

Note that in fact condition (\P) does not depend on the metric on $\unii$, since all metrics on a compact space are uniformly equivalent.

\begin{pf}
First, by Proposition~\ref{Ppicewiselins}, we may assume that all maps $v_n^m$ are piece-wise linear.
Next, by an easy induction, we can ``convert" $\wek v$ to a sequence in $\lcat \fI$.
It is clear that (\P) is preserved under the equivalence of sequences,
therefore the ``corrected" sequence still satisfies (\P).
Now it is obvious that (\P) is equivalent to condition (A) of the \fra\ sequence, since the map $f$ can be approximated by a piece-wise linear quotient map which in turn can be made 1-Lipschitz by multiplying the metric of $\unii$ by a large enough constant.

Thus, if $\wek v$ satisfies (\P) then it is equivalent to a \fra\ sequence which is uniquely determined, showing that $\liminv \wek v \iso P$.
Finally, if $\liminv \wek v \iso P$, then $\wek v$ is equivalent to $\wek u$, therefore it satisfies (\P).
\end{pf}

Lemma~\ref{Legrdin} allows us to work in the monoidal category of quotient maps of the unit interval, endowed with the standard metric.
This is formally not a metric-enriched category, because the composition operator is not 1-Lipschitz, but in practice this does not cause any trouble.

\begin{lm}\label{Lmebrnre}
Every non-degenerate subcontinuum of $P$ is homeomorphic to $P$.
\end{lm}

\begin{pf}
Let $K$ be a subcontinuum of $P$ and let $K_n = \img {u_n}K$, where $\map {u_n} P \unii$ is the canonical $n$-th projection.
From some point on, $K_n$ is a non-degenerate interval.
Without loss of generality, we may assume that this is the case for all $\ntr$.
Given $n<m$, let $\map {v_n^m}{K_m}{K_n}$ be the restriction of $u_n^m$.
Then $\wek v = \sett{v_n^m}{n<m<\nat}$ is a sequence in $\fI$, an inverse sequence of piece-wise linear quotient maps of closed intervals.
Furthermore, $K = \liminv \wek v$.
It suffices to check that $\wek v$ is a \fra\ sequence in $\fI$.

Fix $\ntr$, $\eps > 0$ and fix a piece-wise linear quotient map $\map f \unii F_n$.
Assume $F_n = [a,b]$, where $0 \loe a < b \loe 1$.
Composing $f$ with a suitable quotient map, we may assume that $f(0) = a$ and $f(1) = b$.
Extend $f$ to a piece-wise linear map $\map {f'}{[-1,2]}\unii$ in such a way that $\img {f'}{[-1,0]} = [0,a]$ and $\img {f'}{[1,2]} = [b,1]$.
We can treat $[-1,2]$ as the unit interval with multiplied metric.
Thus, using the fact that $\wek u$ is \fra, we find $m > n$ and a piece-wise linear quotient map $\map g \unii {[-1,2]}$ such that $\rho(f' \cmp g, u_n^m) < \eps$.
Note that $\img g {F_m}$ is $\eps$-close to $[a,b]$, therefore we 
can ``correct" $g$ so that $\img g{F_m} = [a,b]$, replacing $\eps$ by $3 \eps$.
Finally, $g \rest F_m$ witnesses that $\wek v$ satisfies condition (A) of the definition of a \fra\ sequence.
\end{pf}

\begin{lm}
$P$ is hereditarily indecomposable.
\end{lm}

\begin{pf}
In view of Lemma~\ref{Lmebrnre}, it suffices to show that $P$ is indecomposable.
For this aim, suppose $P = A \cup B$, where $A, B$ are proper subcontinua of $P$.
Let $A_n = \img {u_n}A$, $B_n = \img {u_n}B$, where as before, $u_n$ is the canonical $n$th projection.
Fix $\ntr$ such that both $A_n$ and $B_n$ are non-degenerate proper intervals.
Without loss of generality, we may assume that $A_n = [0,a]$, $B_n = [b,1]$, where $0 < b \loe a < 1$.
Let $\map f \unii \unii$ be the tent map, that is, $f(t) = 2t$ for $t \in [0,\frac 12]$ and $f(t) = \frac 12 - 2t$ for $t \in [\frac 12, 1]$.
Then $\inv f {A_n} = [0,s] \cup [t,1]$, where $s = \frac a2$ and $t = 1 - \frac a2 > s$.
Furthermore, $f$ is $2$-Lipschitz with respect to the standard metric, therefore multiplying the metric in the domain of $f$ by $2$ we obtain a 1-Lipschitz piece-wise linear quotient map.
Fix a small enough $\eps>0$.
Using property (A) of the \fra\ sequence, we find $m > n$ and a quotient map $\map g \unii \unii$ such that $\rho(f \cmp g, u_n^m) < \eps$.

Notice that $J = \img g{A_m}$ is an interval, therefore if $\eps$ is small enough then either $J \cap [0,s] = \emptyset$ or $J \cap [t,1] = \emptyset$.
This means that either $[0,s] \subs \img g{B_m}$ or $[t,1] \subs \img g{B_m}$.
In particular, there is $r \in B_m$ such that $g(r) \in \dn 01$.
On the other hand, $d(f(g(r)), u_n^m(r)) = d(0, u_n^m(r)) < \eps$, where $d$ is the metric in the $n$th interval of the sequence $\wek u$. Notice that $u_n^m(r) \in \img {u_n^m}{B_m} = B_n$.
Thus, if $\eps < d(0,b)$ then we get a contradiction.
\end{pf}

The two lemmas above together with Bing's uniqueness result \cite{Bing} give

\begin{wn}
$P$ is the pseudo-arc.
\end{wn}

Applying Theorem~\ref{Thunivweek}, we obtain another proof of the result of Mioduszewski~\cite{Miod2}:

\begin{tw}
Every chainable continuum is a continuous image of the pseudo-arc.
\end{tw}

Almost homogeneity can be easily strengthened, obtaining the result of Irwin \& Solecki \cite{IrwSol} which in turn improves a result of Lewis (sketched after Thm. 4.2 in~\cite{Lewis}):

\begin{tw}\label{Tthgds}
Let $K$ be a chainable continuum with some fixed metric and let $\map {p,q} P K$ be quotient maps.
Then for each $\eps > 0$ there exists a homeomorphism $\map h P P$ such that $\rho(q \cmp h, p) < \eps$.
\end{tw}

\begin{pf}
Using the fact that $K$ is the inverse limit of unit intervals, there is a quotient map $\map f K \unii$ such that all $f$-fibers have diameter $< \eps$.
A standard compactness argument shows that $f$ satisfies the following condition:
\begin{equation}
(\forall \; s,t \in K) \;\; |f(s) - f(t)| < \delta \implies d(s,t) < \eps.
\tag{$\star$}\label{Eqsterrv}
\end{equation}
Now let $p' = f \cmp p$ and $q' = f \cmp q$.
By Lemma~\ref{Lewfkt}, both $p'$, $q'$ come from approximate arrows, therefore by Theorem~\ref{TwAlmHomgy}, there is a homeomorphism $\map h P P$ such that $\rho(q' \cmp h, p') < \delta$.
We have the following diagram, in which the upper triangle is $\delta$-commutative and the side-triangles are commutative.
$$\xymatrix{
P \ar[rr]^h \ar[dd]_p \ar[ddr]^{p'} & & P \ar[dd]^q \ar[ddl]_{q'} \\
& & \\
K \ar[r]_f & \unii & K \ar[l]^f
}$$

Finally, condition (\ref{Eqsterrv}) gives $\rho(q \cmp h, p) < \eps$.
\end{pf}

As one can guess, the property in Theorem~\ref{Tthgds} characterizes the pseudo-arc among chainable continua.
In fact, this has already been proved by Irwin \& Solecki~\cite{IrwSol}.

\begin{tw}
A chainable continuum $K$ is homeomorphic to $P$ if and only if it satisfies the following condition:
\begin{enumerate}
	\item[{\rm(P)}] Given $\eps > 0$, given quotient maps $\map q K \unii$, $\map f \unii \unii$, there exists a quotient map $\map g K \unii$ such that
$\rho(f \cmp g, q) < \eps$.
\define{}{(P)}
\end{enumerate}
\end{tw}

\begin{pf}
By Theorem~\ref{Tthgds}, $P$ satisfies condition (P).
Now suppose that $K$ satisfies (P) and choose a sequence $\wek v$ in $\fI$ whose inverse limit is $K$.
We shall check that $\wek v$ is a \fra\ sequence.
In fact, only condition (A) requires a proof.

Fix $\ntr$, $\eps > 0$ and fix a piece-wise linear quotient map $\map f \unii \unii$.
Let, as usual, $\map{v_n}K \unii$ be the $n$th canonical projection.
Using (P), we find a quotient map $\map p K \unii$ such that $\rho(f \cmp p, v_n) < \eps/2$.
By Lemma~\ref{Lewfkt}, there exist $m > n$ and a piece-wise linear quotient map $\map g \unii \unii$ such that $\rho(g \cmp v_m, p) < \eps/2$.
Thus we get
$$\rho(f \cmp g \cmp v_m, v_n^m \cmp v_m) \loe \rho(f \cmp g \cmp v_m, f \cmp p) + \rho(f \cmp p, v_n) < \eps$$
and so $\rho(f \cmp g, v_n^m) < \eps$, because $v_m$ is a quotient map.
This shows (A) and completes the proof.
\end{pf}

As noticed at the beginning of the proof above, condition (P) easily follows from Theorem~\ref{Tthgds}.
A direct proof of the converse implication would require the approximate back-and-forth argument, which is hidden in the proof of Lemma~\ref{Labnf} above.

\begin{uwgi}
As indicated in \cite{IrwSol}, one can try to prove the existence, properties and uniqueness of the so-called \emph{pseudo-circle}, thus extending the work of Rogers~\cite{Rogers}.
It seems that this can be done using piece-wise linear non-zero degree self-maps of the circle, for which the uniformization theorem was proved by Rogers~\cite{Rogers}.
\end{uwgi}



\printindex

\end{document}